\documentclass[journal]{IEEEtran}

\usepackage{amsmath}
\usepackage{amsfonts}
\usepackage{gensymb}
\usepackage{graphicx}
\usepackage{float}
\usepackage[small,bf]{caption}
\usepackage{url}
\usepackage{color}
\usepackage{epstopdf}
\usepackage{cuted}
\usepackage{tikz}
\usepackage{tabularx}
\usepackage{tabulary}
\usepackage{todonotes}
\usepackage{lipsum}
\usepackage{enumitem}
\usepackage{multirow}
\usepackage{url}
\usepackage{empheq}
\usepackage[most]{tcolorbox}
\usepackage[export]{adjustbox}
\usepackage{amssymb}

\graphicspath{{figures/},{pdf/},{eps/},{png/}}

\renewcommand{\Re}{\operatorname{\textbf{Re}}}
\renewcommand{\Im}{\operatorname{\textbf{Im}}}

\newcommand{\q}{\,}
\newcommand{\qq}{\,\,}

\definecolor{mnavy}{RGB}{12,44,86}
\definecolor{mblue}{RGB}{0,40,120}

\tcbset{colback=white, colframe=mblue, sharp corners, left=-9.25pt, right=-0.75pt, top=-5pt, bottom=0pt}

\begin{document}

\title{Optimal Voltage Phasor Regulation for Switching Actions in Unbalanced Distribution Systems}

\author{\IEEEauthorblockN{Michael D. Sankur, Roel Dobbe, Alexandra von Meier, Emma M. Stewart, and Daniel B. Arnold}
}

\maketitle

\let\thefootnote\relax\footnote{Michael D. Sankur and Daniel B. Arnold are with Lawrence Berkeley National Laboratory. Roel Dobbe and Alexandra Von Meier are with the Electrical Engineering Department at U.C. Berkeley. Emma Stewart is with the Lawrence Livermore National Laboratory. This work was supported in part by the U.S. Department of Energy ARPA-E program (DE-AR0000340), and the Department of Energy Office of Energy Efficiency and Renewable Energy under Contract No. DE-AC02-05CH11231.}

\begin{abstract}

The proliferation of phasor measurement units (PMUs) into electric power distribution grids presents new opportunities for utility operators to manage distribution systems more effectively.
One potential application of PMU measurements is to facilitate distribution grid re-configuration.
Given the increasing amount of Distributed Energy Resource (DER) penetration into distribution grids, in this work we formulate an Optimal Power Flow (OPF) approach that manages DER power injections to minimize the voltage phasor difference between two nodes on a distribution network to enable efficient network reconfiguration.
In order to accomplish this, we develop a linear model that relates voltage phase angles to real and reactive power flows in unbalanced distribution systems.
Used in conjunction with existing linearizations relating voltage magnitude differences to power flows, we formulate an OPF capable of minimizing voltage phasor differences across different points in the network.
In simulations, we explore the use of the developed approach to minimize the phasor difference across switches to be opened or closed,
thereby providing an opportunity to automate and increase the speed of reconfigurations in unbalanced distribution grids.

\end{abstract}



\section{Introduction}
\label{sec:introduction}

The proliferation of new types of sensors into the electric power distribution system is providing deeper insights into grid operation and is driving innovation around new paradigms for system management. Among the many new sensing devices being deployed in the distribution grid, distribution phasor measurement units (PMUs) provide a suite of new functionality that could serve to better inform the process of managing Distributed Energy Resources (DER). Distribution PMUs provide access to the magnitude and angle of voltage and current phasors. These devices are now becoming more commonplace and are either manifesting as standalone units \cite{vonMeier2017precision} or are being incorporated into other system components \cite{synchrophasor2010}. 

The expansion of PMU technology from transmission into distribution systems indicates that sufficient infrastructure may be in place in future grids to support control activities that make decisions based on feeder voltage phasor measurements. In fact, a small, but growing, number of control applications that utilize phase angle measurements have started to appear in literature. The work of \cite{ochoa2010angle} proposed the use of synchronized voltage phase angle measurements to curtail over-generation of renewables. Additionally, the authors of \cite{wang2013pmu} considered voltage angle thresholds as criteria to connect renewable generation. Both works refer to this control activity as ``Angle Constrained Active Management'', or ACAM.

One particular emerging application for which distribution PMU data might be of use is to enable fast and safe switching of circuit elements. The ability to island/reconnect microgrids and reconfigure distribution feeders are seen as two important applications of future grids \cite{grid2015,quad2015,heydt2010next}. 
In order to close a switch, the voltage magnitude and voltage angle at both sides should be sufficiently close in order to prevent arcing and transient currents.
As such, to facilitate network reconfigurations, distribution system operators (DSOs) typically employ backup power injection at reconfiguration locations to minimize voltage phasor difference across switches.
DSOs either schedule such actions in advance, in the case where a reconfiguration is planned, and send mobile generators with back up power to the correct locations.
In the event of an unplanned outage, it may take considerable time to deploy engineering staff with backup power to the necessary switching locations.

Distributed energy resources (DER) present an opportunity to facilitate network reconfiguration at timescales much faster than with current practices. With sufficient amount and proper location of DERs it is possible to control the voltage phasor at strategic points in the network, thereby alleviating the need for ad-hoc deployment of generation assets for switching. The ability to switch elements in and out of a given system with the aid of controlled DERs could allow for faster restoration of electrical services to critical loads following a disaster, or allow for damaged components to be isolated for repair or replacement. 

In literature, there is no shortage of new strategies that consider the use of DER to operate distribution feeders in a more sophisticated manner. Many of the approaches formulate the decision making process as an optimization program, often referred to as an Optimal Power Flow (OPF) problem. An OPF framework allows for proper modeling of the network topology, impedances and control equipment, the incorporation of safety constraints, and the formulation of various objective functions that can reflect important operating objectives such as loss minimization and cost of generation or control action. However, there is sparse literature on strategies that aim to \emph{directly control the voltage phasor in unbalanced systems in OPF formulations}. The work of \cite{zhang2017distributed} proposes a distributed control framework to enable DER to track single-phase AC optimal power flow solutions using the Alternating Direction Method of Multipliers (ADMM). The authors of \cite{xu2017multi} consider a multi-timescale stochastic volt/var control method capable of controlling legacy voltage regulation systems as well as DER. Some of the authors of the present work consider optimal governance of DER in a model-free setting \cite{arnold2017model}. Additional recent strategies for DER control are cataloged in \cite{antoniadou2017distributed}.

Due to the nonlinear nature of power flow equations, many OPFs are formulated as quadratically constrained quadratic programs (QCQPs). A popular method for analyzing such OPFs is relaxation via semidefinite programming (SDP) \cite{dall2012optimization,dall2013distributed}. It is well documented that relaxation of OPF problems via SDP often fails to achieve a rank-one solution \cite{lesieutre2011examining,louca2014nondegeneracy,madani2015convex}.

As an example, in the work of \cite{louca2014nondegeneracy} too many binding constraints will preclude convergence to a rank-one solution. The authors of \cite{madani2015convex} explored the extension of SDP to weakly meshed networks. Their technique was able to achieve a rank-one solution only after incorporating significant penalties on reactive power dispatch, effectively limiting the feasible region of control. Finally, in \cite{dall2013distributed}, the authors faced difficulty in obtaining a rank-one solution for certain network configurations.

As the inability of relaxations via SDPs to achieve a rank one solution limits the practicality of these approaches, it is necessary to consider alternative approaches for solving OPF problems.  One such alternative is the creation of linear approximations for power flow that are sufficiently accurate for control purposes, and that can be incorporated into convex OPF formulations.




\subsection*{Contributions}
To enable a control strategy that can regulate voltage phasors, in this work we extend a linearized model of three phase power flow to include a mapping of voltage phase angle differences into network real and reactive power flows. We do this, as to our knowledge, OPF approaches controlling voltage magnitude, and active and reactive power flows cannot always minimize voltage phasor difference across an open switch. Despite known difficulties associated with SDPs, we had first attempted to formulate this activity as an semidefinite-relaxed OPF, but were unsuccessful. This motivated the extension of the linear mappings of \cite{gan2014convex,robbins2016optimal,arnold2015optimal}, which relate voltage magnitude to active and reactive power flows, to consider the entire voltage phasor. Our contribution is an extension of these results to include a linear relation between complex line power flow and voltage angle difference, in Section \ref{sbsec:angle}. Additionally, we formulate an OPF capable of managing voltage phasors, as opposed to only voltage magnitudes.

The specific activity studied herein is an OPF formulation that minimizes the voltage phasor difference across an open switch in a distribution system while simultaneously regulating feeder voltage magnitudes to within acceptable limits. In the event that one of the phasors is uncontrolled (i.e. a reference signal), then this activity can be thought of as a voltage phasor tracking problem. In driving the voltage phasor difference across a circuit element to 0, we ensure that when the switch is closed only small amounts of power will flow across this element. In this manner, the switch can be closed with less arcing and instantaneous power flow surges.

The analysis in this paper has two major elements. A derivation of a linearized model of unbalanced power flow that maps voltage phasor differences into active and reactive power flows is presented in Section \ref{sec:modelderivation}.
In Section \ref{sec:montecarlo}, a numerical analysis is provided that outlines errors introduced by the linearized model. Simulation results of an OPF that incorporates the linearized power flow model to track a voltage phasor reference at a specific point in the network and regulate system voltage magnitudes are then presented in Section \ref{sec:phasortracking}.

\section{Model Derivation}
\label{sec:modelderivation}

In this section, we derive the linearized unbalance power flow model. Please see the nomenclature table for variable definitions used in subsequent analysis.

\subsection{Preliminaries}
\label{sbsec:preliminaries}

\renewcommand{\arraystretch}{1.15}
\setlength{\textfloatsep}{10pt}
\begin{table}[t]
\centering
\caption*{NOMENCLATURE}
\begin{tabulary}{\linewidth}{LLL}
\hline
\hline
$\mathcal{P}_{n}$ & & Set of phases that exist at node $n$ \\
$\mathcal{P}_{mn}$ & & Set of phases that exist on line $(m,n)$ \\
$V_{n}^{\phi}$ & & Voltage phasor on phase $\phi$ at node $n$ \\
$\mathbf{V}_{n}$ & & Vector of voltage phasors at node $n$ \\
$E_{n}^{\phi}$ & & Squared magnitude of voltage on phase $\phi$ at node $n$ \\
$\mathbf{E}_{n}$ & & Vector of squared magnitudes of voltage at node $n$ \\
$\theta_{n}^{\phi}$ & & Angle of voltage phasor on phase $\phi$ at node $n$\\
$\Theta_{n}$ & & Vector of voltage phasor angles at node $n$ \\
$Z_{mn}^{\phi \psi}$ & & Impedance of line $(m,n)$ between phases ($\phi, \psi$) \\
$\mathbf{Z}_{mn}$ & & Impedance matrix of line $(m,n)$ \\
$I_{mn}^{\phi}$ & & Current phasor on phase $\phi$ on line $(m,n)$  \\
$\mathbf{I}_{mn}$ & & Vector of current phasors on line $(m,n)$ \\
$i_{n}^{\phi}$ & & Node current on phase $\phi$ at node $n$ \\
$\mathbf{i}_{n}$ & & Vector of node currents at node $n$ \\
$S_{mn}^{\phi}$ & & Phasor of complex power entering node $n$ on phase $\phi$ on line $(m,n)$ \\
$\mathbf{S}_{mn}$ & & Vector of complex power phasors entering node $n$ on line $(m,n)$ \\
$s_{n}^{\phi}$ & & Complex nodal power phasor on phase $\phi$ at node $n$ \\
$\mathbf{s}_{n}$ & & Vector of complex nodal power phasors at node $n$ \\
$w_{n}^{\phi}$ & & Controllable complex power dispatch on phase $\phi$ at node $n$ \\
$\left( \cdot \right)^{*}$ & & Complex conjugate \\
\hline
\hline
\end{tabulary}
\end{table}

Let $\mathcal{T} = (\mathcal{N}, \mathcal{E})$ denote a graph representing an unbalanced distribution feeder, where $\mathcal{N}$ is the set of nodes of the feeder and $\mathcal{E}$ is the set of line segments. Nodes are indexed by $m$ and $n$, with $m$, $n \in \mathcal{N}$.
Let $\mathcal{N} \triangleq \{ \infty, 0, 1, \ldots \vert \mathcal{N} \vert \}$, where node 0 denotes the substation (feeder head). Immediately upstream of node 0 is an additional
node used to represent the transmission system, indexed by $\infty$.  We treat node $\infty$ as an infinite bus, decoupling interactions in the downstream distribution system from the rest of the grid. While the substation voltage may evolve over time, we assume this evolution takes place independently of DER control actions in $\mathcal{T}$.

Each node and line segment in $\mathcal{T}$ can have up to three phases, labeled $a$, $b$, and $c$. Phases are referred to by $\phi \in \{ a, b, c \}$ and $\psi \in \{ a, b, c \}$. We define $\mathcal{P}_{m}$ and $\mathcal{P}_{n}$ as the set of phases at nodes $m$ and $n$, respectively, and $\mathcal{P}_{mn}$ as set of phases of line segment $(m,n)$. If phase $\phi$ is present at node $m$, then at least one line connected to $m$ must contain phase $\phi$.
If line $(m,n)$ exists, its phases are a subset of the phases present at both node $m$ and node $n$, such that $(m,n) \in \mathcal{E} \Rightarrow \mathcal{P}_{mn} \subseteq \mathcal{P}_{m} \cap \mathcal{P}_{n} $.

The current/voltage relationship for a three phase line $(m,n)$ between adjacent nodes $m$ and $n$ is captured by Kirchhoff's Voltage Law (KVL) in its full \eqref{eqn:KVLmn01}, and vector form \eqref{eqn:KVLmn02}:
\begin{gather}
\begin{bmatrix}
V_{m}^{a} \\ V_{m}^{b} \\ V_{m}^{c}
\end{bmatrix}
=
\begin{bmatrix}
V_{n}^{a} \\ V_{n}^{b} \\ V_{n}^{c}
\end{bmatrix}
+
\begin{bmatrix}
Z^{aa}_{mn} & Z^{ab}_{mn} & Z^{ac}_{mn} \\
Z^{ba}_{mn} & Z^{bb}_{mn} & Z^{bc}_{mn} \\
Z^{ca}_{mn} & Z^{cb}_{mn} & Z^{cc}_{mn}
\end{bmatrix}
\begin{bmatrix}
I_{mn}^{a} \\ I_{mn}^{b} \\ I_{mn}^{c}
\end{bmatrix}
\q ,
\label{eqn:KVLmn01}
\\
\mathbf{V}_{m} = \mathbf{V}_{n} + \mathbf{Z}_{mn} \mathbf{I}_{mn} \q .
\label{eqn:KVLmn02}
\end{gather}
Here, $Z^{\phi \psi}_{mn} = r^{\phi \psi}_{mn} + jx^{\phi \psi}_{mn}$ denotes the complex impedance of line $(m,n)$ across phases $\phi$ and $\psi$. We have presented \eqref{eqn:KVLmn01} and \eqref{eqn:KVLmn02} where $\mathcal{P}_{mn} = \{ a, b, c\}$. For lines with fewer than three phases ($\left| \mathcal{P}_{mn} \right| \le 2$) \eqref{eqn:KVLmn02} becomes:
\begin{equation}
\left[
\mathbf{V}_{m}
=
\mathbf{V}_{n}
+
\mathbf{Z}_{mn}
\mathbf{I}_{mn}
\right]_{\mathcal{P}_{mn}}
\q ,
\label{eqn:KVLmn03}
\end{equation}

\noindent by indexing by the set of line phases $\mathcal{P}_{mn}$, where the rows associated with phases $\psi \notin \mathcal{P}_{mn}$ of \eqref{eqn:KVLmn02} are removed, as are the appropriate columns of $\mathbf{Z}_{mn}$. To give two examples, if $\mathcal{P}_{mn} = \left\{ a \right\}$, then \eqref{eqn:KVLmn03} is $\left[ \mathbf{V}_{m} \right]_{ \left\{ a \right\} } \equiv V_{m}^{a} = V_{n}^{a} + Z_{mn}^{aa} I_{mn}^{a}$, and if $\mathcal{P}_{mn} = \left\{ a , c \right\}$ then \eqref{eqn:KVLmn03} is:
\begin{equation*}
\left[ \mathbf{V}_{m} \right]_{ \left\{ a, c \right\} }
\equiv
\begin{bmatrix}
V_{m}^{a} \\ V_{m}^{c}
\end{bmatrix}
=
\begin{bmatrix}
V_{n}^{a} \\ V_{n}^{c}
\end{bmatrix}
+
\begin{bmatrix}
Z^{aa}_{mn} & Z^{ac}_{mn} \\
Z^{ca}_{mn} & Z^{cc}_{mn}
\end{bmatrix}
\begin{bmatrix}
I_{mn}^{a} \\ I_{mn}^{c}
\end{bmatrix}
\q .
\end{equation*}

Kirchoff's Current law at node $m$ is given in its full \eqref{eqn:KCLm01} and vector \eqref{eqn:KCLm02} forms:
\begin{gather}
\sum_{l:(l,m) \in \mathcal{E}}
\begin{bmatrix}
I_{lm}^{a} \\ I_{lm}^{b} \\ I_{lm}^{c}
\end{bmatrix}
=
\begin{bmatrix}
i_{m}^{a} \\ i_{m}^{b} \\ i_{m}^{c}
\end{bmatrix}
+
\sum_{n:(m,n) \in \mathcal{E}}
\begin{bmatrix}
I_{mn}^{a} \\ I_{mn}^{b} \\ I_{mn}^{c}
\end{bmatrix}
\q ,
\label{eqn:KCLm01}
\\
\sum_{l:(l,m) \in \mathcal{E}} \mathbf{I}_{lm}
=
\mathbf{i}_{m}
+
\sum_{n:(m,n) \in \mathcal{E}} \mathbf{I}_{mn} 
\q .
\label{eqn:KCLm02}
\end{gather}


\noindent We assume a complex load, $s_{n}^{\phi} \qq \forall \phi \in \mathcal{P}_{n}, \qq \forall n \in \mathcal{N} \setminus \infty$, is served on all existing phases at each node except the transmission line, defined as:
\begin{equation}
s_{n}^{\phi} \left( V_{n}^{\phi} \right) = \left( \beta_{S,n}^{\phi} + \beta_{Z,n}^{\phi} \left| V_{n}^{\phi} \right|^{2} \right) d_{n}^{\phi} + w_{n}^{\phi} - j c_{n}^{\phi} \q ,
\label{eqn:sV}
\end{equation}

\noindent where $d_{n}^{\phi}$ is the complex demand, with constant power and impedance terms $\beta_{S,n}^{\phi} + \beta_{Z,n}^{\phi} = 1$, $w_{n}^{\phi} = u_{n}^{\phi} + j v_{n}^{\phi}$ represents complex power available for control (e.g. DER), and $c_{n}^{\phi}$ denotes capacitance, all for $\phi \in \mathcal{P}_{n}$. Note, if $\phi \notin \mathcal{P}_{n}$ (i.e. phase $\phi$ does not exist at node $n$), we define $V_{n}^{\phi} = i_{n}^{\phi} = s_{n}^{\phi} = d_{n}^{\phi} = w_{n}^{\phi} = c_{n}^{\phi} = 0$. If $\phi \notin \mathcal{P}_{mn}$ (i.e. phase $\phi$ does not exist on line segment $(m,n)$), we define $I_{mn}^{\phi} = 0$.

Throughout this work, we use the symbol $\circ$ to represent the Hadamard Product (HP) of two matrices of the same dimension, also known as the element-wise product, which can be written as:
\begin{equation*}
C = A \circ B = B \circ A \Rightarrow C_{ij} = A_{ij} B_{ij} = B_{ij} A_{ij} \q .
\label{eqn:hadamard}
\end{equation*}


\subsection{Power and Losses}
\label{sbsec:power}

We now derive complex power and loss terms at a node $m \in \mathcal{N}$.
This analysis, and the derivation of \ref{sbsec:magnitude}, we do not claim as novel contributions (see \cite{gan2014convex,robbins2016optimal}).
Full derivation of these results are necessary to support one of the main contributions of this work, which is presented in \ref{sbsec:angle}.
To start, we take the Hadamard Product of $\mathbf{V}_{m}$ and the complex conjugate (non-transposed) of \eqref{eqn:KCLm02}:
\begin{equation}
\sum_{l:(l,m) \in \mathcal{E}} \mathbf{V}_{m} \circ \mathbf{I}_{lm}^{*}
=
\mathbf{V}_{m} \circ \mathbf{i}_{m}^{*}
+
\sum_{n:(m,n) \in \mathcal{E}} \mathbf{V}_{m} \circ \mathbf{I}_{mn}^{*} \q .
\label{eqn:power01}
\end{equation}

\noindent The $\mathbf{V}_{m}$ term inside the summation on the RHS is substituted using \eqref{eqn:KVLmn02}:
\begin{equation}
\begin{aligned}
\sum_{l:(l,m) \in \mathcal{E}} \mathbf{V}_{m} & \circ \mathbf{I}_{lm}^{*}
=
\mathbf{V}_{m} \circ \mathbf{i}_{m}^{*} + \ldots \\
& \sum_{n:(m,n) \in \mathcal{E}} \mathbf{V}_{n} \circ \mathbf{I}_{mn}^{*}
+
\left( \mathbf{Z}_{mn} \mathbf{I}_{mn} \right) \circ \mathbf{I}_{mn}^{*}
\q .
\end{aligned}
\label{eqn:power02}
\end{equation}

\noindent Here, we define the complex power phasor on phase $\phi$ entering node $n$ on line $(m,n)$ as $S_{mn}^{\phi} = V_{n}^{\phi} \left( I_{mn}^{\phi} \right)^{*}$, and the $3 \times 1$ vector of complex power phasors entering node $n$ on line $(m,n)$ as $\mathbf{S}_{mn} = \left[ S_{mn}^{a}, \qq S_{mn}^{b}, \qq S_{mn}^{c} \right]^{T} = \mathbf{V}_{n} \circ \mathbf{I}_{mn}^{*}$.
The complex load at node $m$ on phase $\phi$ is defined as $s_{m}^{\phi} = V_{m}^{\phi} \left( i_{m}^{\phi} \right)^{*}$, and the $3 \times 1$ vector of complex load phasors at node $m$ is $\mathbf{s}_{m} = \left[ s_{m}^{a}, \qq s_{m}^{b}, \qq s_{m}^{c} \right]^{T} =\mathbf{V}_{m} \circ \mathbf{i}_{m}^{*}$.
We now rewrite \eqref{eqn:power02}:
\begin{equation}
\sum_{l:(l,m) \in \mathcal{E}} \mathbf{S}_{lm} = \mathbf{s}_{m} + \sum_{n:(m,n) \in \mathcal{E}} \mathbf{S}_{mn} + \mathbf{L}_{mn} \q .
\label{eqn:power03}
\end{equation}

The term $\mathbf{L}_{mn} \in \mathbf{C}^{3 \times 1}$ represents nonlinear losses on the line.
As in \cite{gan2014convex, robbins2016optimal, baran1989optimal}, we assume that losses are negligible compared to line flows, so that $\left| L_{mn}^{\phi} \right| \ll \left| S_{mn}^{\phi} \right| \qq \forall (m,n) \in \mathcal{E}$.
Thus, we neglect line losses, linearizing \eqref{eqn:power03} into:
\begin{equation}
\sum_{l:(l,m) \in \mathcal{E}} \mathbf{S}_{lm} \approx \mathbf{s}_{m} + \sum_{n:(m,n) \in \mathcal{E}} \mathbf{S}_{mn} \q .
\label{eqn:power04}
\end{equation}

\subsection{Voltage Magnitude Equations}
\label{sbsec:magnitude}

In this section, we derive a relation between squared voltage magnitudes and complex multiphase power for unbalanced systems.
The reader should note that here we present the derivation for a line with three phases, where $\mathcal{P}_{mn} = \{ a, b, c \}$. For lines with less than three phases ($\left| \mathcal{P}_{mn} \right| \le 2$), \eqref{eqn:mag01} - \eqref{eqn:MmnNmn} should be indexed by $\mathcal{P}_{mn}$, as in \eqref{eqn:KVLmn03}.

To start, we consider a line $(m,n) \in \mathcal{E}$, and take the Hadamard Product of \eqref{eqn:KVLmn02} and its (non-transposed) complex conjugate:
\begin{equation}
	\mathbf{V}_{m} \circ \mathbf{V}_{m}^{*} = \left( \mathbf{V}_{n} + \mathbf{Z}_{mn} \mathbf{I}_{mn} \right) \circ \left( \mathbf{V}_{n} + \mathbf{Z}_{mn} \mathbf{I}_{mn} \right)^{*} \q .
    \label{eqn:mag01}
\end{equation}

\noindent This can be rewritten by distributing the terms on the RHS:
\begin{equation}
\begin{aligned}
	\mathbf{V}_{m} & \circ \mathbf{V}_{m}^{*} = \mathbf{V}_{n} \circ \mathbf{V}_{n}^{*} + \mathbf{V}_{n} \circ  \left( \mathbf{Z}_{mn} \mathbf{I}_{mn} \right)^{*} + \ldots \\
    & \left( \mathbf{Z}_{mn} \mathbf{I}_{mn} \right) \circ \mathbf{V}_{n}^{*} + \left( \mathbf{Z}_{mn} \mathbf{I}_{mn} \right) \circ \left( \mathbf{Z}_{mn} \mathbf{I}_{mn} \right)^{*}  \q .
\end{aligned}
\label{eqn:mag02}
\end{equation}

\noindent Here we define the real scalar $E_{n}^{\phi} = \left| V_{n}^{\phi} \right|^{2} = V_{n}^{\phi} (V_{n}^{\phi})^{*}$, the $3 \times 1$  real vector $\mathbf{E}_{n} = \left[ E_{n}^{a} , \qq E_{n}^{b} , \qq E_{n}^{c} \right]^{T} = \mathbf{V}_{n} \circ \mathbf{V}_{n}^{*}$, and the $3 \times 1$ real vector $\mathbf{H}_{mn} = \left( \mathbf{Z}_{mn} \mathbf{I}_{mn} \right) \circ \left( \mathbf{Z}_{mn} \mathbf{I}_{mn} \right)^{*} = \left( \mathbf{V}_{m} - \mathbf{V}_{n} \right) \circ \left( \mathbf{V}_{m} - \mathbf{V}_{n} \right)^{*}$. With these definitions, we also take advantage of the commutative property of the HP and group the second and third terms of the RHS of \eqref{eqn:mag02} inside the real operator:
\begin{equation}
	\mathbf{E}_{m} = \mathbf{E}_{n} +
	2 \Re \left\{ \left( \mathbf{Z}_{mn} \mathbf{I}_{mn} \right)^{*} \circ \mathbf{V}_{n} \right\} + \mathbf{H}_{mn} \q .
	\label{eqn:mag03}
\end{equation}



At this point, we focus on the terms inside the real operator for clarity of presentation, and rewrite them as: 
\begin{equation}
\begin{aligned}
	& \left( \mathbf{Z}_{mn} \mathbf{I}_{mn} \right)^{*} \circ \mathbf{V}_{n} = \ldots \\
    & \qq
    \begin{bmatrix}
    	V_{n}^{a} \left( Z_{mn}^{aa} I_{mn}^{a} + Z_{mn}^{ab} I_{mn}^{b} + Z_{mn}^{ac} I_{mn}^{c} \right)^{*} \\
        V_{n}^{b} \left( Z_{mn}^{ba} I_{mn}^{a} + Z_{mn}^{bb} I_{mn}^{b} + Z_{mn}^{bc} I_{mn}^{c} \right)^{*} \\
        V_{n}^{c} \left( Z_{mn}^{ca} I_{mn}^{a} + Z_{mn}^{cb} I_{mn}^{b} + Z_{mn}^{cc} I_{mn}^{c} \right)^{*}
    \end{bmatrix} \q .
\end{aligned}
\label{eqn:mag04}
\end{equation}

\noindent With the definition of complex current on a line, $ I_{mn}^{\phi} = \left( S_{mn}^{\phi} / V_{n}^{\phi} \right)^{*}$, and defining the term $\gamma_{n}^{\phi \psi} = V_{n}^{\phi} / V_{n}^{\psi}$, we rewrite \eqref{eqn:mag04} as:
\begin{equation}
\begin{aligned}
	& \left( \mathbf{Z}_{mn} \mathbf{I}_{mn} \right)^{*} \circ \mathbf{V}_{n} = \ldots
    \\
    & \qq
    \begin{bmatrix}
    	(Z_{mn}^{aa})^{*} S_{mn}^{a} + \gamma_{n}^{ab} (Z_{mn}^{ab})^{*} S_{mn}^{b} + \gamma_{n}^{ac} (Z_{mn}^{ac})^{*} S_{mn}^{c} \\
        \gamma_{n}^{ba} (Z_{mn}^{ba})^{*} S_{mn}^{a} + (Z_{mn}^{bb})^{*} S_{mn}^{b} + \gamma_{n}^{bc} (Z_{mn}^{bc})^{*} S_{mn}^{c} \\
        \gamma_{n}^{ca} (Z_{mn}^{ca})^{*} S_{mn}^{a} + \gamma_{n}^{cb} (Z_{mn}^{cb})^{*} S_{mn}^{b} + (Z_{mn}^{cc})^{*} S_{mn}^{c}
    \end{bmatrix} \q .
\end{aligned}
\label{eqn:mag05}
\end{equation}

\noindent The $3 \times 1$ vector on the RHS of \eqref{eqn:mag05} can be separated into a $3 \times 3$ matrix multiplying the $3 \times 1$ vector $\mathbf{S}_{mn}$:
\begin{equation}
\begin{aligned}
	& \left( \mathbf{Z}_{mn} \mathbf{I}_{mn} \right)^{*} \circ \mathbf{V}_{n} = \ldots
    \\
    & \qq
    \begin{bmatrix}
    	(Z_{mn}^{aa})^{*} & \gamma_{n}^{ab} (Z_{mn}^{ab})^{*} & \gamma_{n}^{ac} (Z_{mn}^{ac})^{*} \\
        \gamma_{n}^{ba} (Z_{mn}^{ba})^{*} & (Z_{mn}^{bb})^{*} & \gamma_{n}^{bc} (Z_{mn}^{bc})^{*}\\
        \gamma_{n}^{ca} (Z_{mn}^{ca})^{*} & \gamma_{n}^{cb} (Z_{mn}^{cb})^{*} & (Z_{mn}^{cc})^{*}
    \end{bmatrix}
    \begin{bmatrix}
    	S_{mn}^{a} \\ S_{mn}^{b} \\ S_{mn}^{c}
    \end{bmatrix} \q .
\end{aligned}
\label{eqn:mag06}
\end{equation}

We use the definition of the HP to factor the $3 \times 3$ matrix into two $3 \times 3$ matrices as in \eqref{eqn:mag07}, where $\Gamma_{n}$ is the $3 \times 3$ matrix to the left of the Hadamard Product symbol $( \circ )$ within the parentheses on the RHS:
\begin{equation}
\begin{aligned}
	& \left( \mathbf{Z}_{mn} \mathbf{I}_{mn} \right)^{*} \circ \mathbf{V}_{n}
    =
    \left( \Gamma_{n} \circ \mathbf{Z}_{mn}^{*} \right) \mathbf{S}_{mn} = \ldots
    \\
    & \qq
    \left(
    \begin{bmatrix}
    	1 & \gamma_{n}^{ab} & \gamma_{n}^{ac} \\
        \gamma_{n}^{ba} & 1 & \gamma_{n}^{bc} \\
        \gamma_{n}^{ca} & \gamma_{n}^{cb} & 1
    \end{bmatrix}
    \circ
    \begin{bmatrix}
    	Z_{mn}^{aa} & Z_{mn}^{ab} & Z_{mn}^{ac} \\
        Z_{mn}^{ba} & Z_{mn}^{bb} & Z_{mn}^{bc} \\
        Z_{mn}^{ca} & Z_{mn}^{cb} & Z_{mn}^{cc}
    \end{bmatrix}^{*}
    \right)
    \begin{bmatrix}
    	S_{mn}^{a} \\ S_{mn}^{b} \\ S_{mn}^{c}
    \end{bmatrix} \q .
\end{aligned}
\label{eqn:mag07}
\end{equation}


Placing \eqref{eqn:mag07} back into \eqref{eqn:mag03} gives:
\begin{equation}
\begin{aligned}
	\mathbf{E}_{m} = \mathbf{E}_{n} + 2 \Re \left\{ \left( \Gamma_{n} \circ \mathbf{Z}_{mn}^{*} \right) \mathbf{S}_{mn} \right\} + \mathbf{H}_{mn} \q .
\end{aligned}
\label{eqn:mag08}
\end{equation}

\noindent Finally, we separate the complex power vector into its active and reactive components, $\mathbf{S}_{mn} = \mathbf{P}_{mn} + j \mathbf{Q}_{mn}$,
and apply the real operator on the RHS to obtain
\begin{equation}
\begin{gathered}
	\mathbf{E}_{m} = \mathbf{E}_{n} + 2 \mathbf{M}_{mn} \mathbf{P}_{mn} - 2 \mathbf{N}_{mn} \mathbf{Q}_{mn} + \mathbf{H}_{mn} \\
    \mathbf{M}_{mn} = \Re \left\{ \Gamma_{n} \circ \mathbf{Z}_{mn}^{*} \right\}, \q
    \mathbf{N}_{mn} = \Im \left\{ \Gamma_{n} \circ \mathbf{Z}_{mn}^{*} \right\} \q .
\end{gathered}
\label{eqn:mag09}
\end{equation}

We have derived equations that govern the relationship between squared voltage magnitudes and complex power flow across line $(m,n)$. This nonlinear and nonconvex system is difficult to directly incorporate into an OPF formulation without the use of convex relaxations. Following the analysis in \cite{gan2014convex}, we apply two approximations. The first is that the higher order term $\mathbf{H}_{mn}$, which is the change in voltage associated with losses, is negligible, implying $\mathbf{H}_{mn} \approx \left[ 0, \qq 0, \qq 0 \right]^{T} \qq \forall (m,n) \in \mathcal{E}$ \cite{gan2014convex}.
The second fixes the ratio of voltages phasors in different network phases as constants.  The impact of this assumption was investigated in \cite{gan2014convex} where it was shown to result in relatively small errors. Effective choices for the values of the voltage ratios are $1 \pm \angle 120^{\circ}$ which accurately approximate typical operating conditions in distribution systems. With these assumptions, $\Gamma_n$  becomes:
\begin{equation}
	\Gamma_{n}
    \approx
    A
    =
    \begin{bmatrix}
    	1 & \alpha & \alpha^{2} \\
        \alpha^{2} & 1 & \alpha \\
        \alpha & \alpha^{2} & 1
    \end{bmatrix}
	\qq \forall n \in \mathcal{N} \q ,
    \label{eqn:gammaalpha}
\end{equation}

\noindent where $\alpha = 1 \angle 120 \degree = \frac{1}{2} ( -1 + j \sqrt{3} )$ and $\alpha^{2} = \alpha^{-1} = \alpha^{*} = 1 \angle 240 \degree = \frac{1}{2} ( -1 - j \sqrt{3} )$. We note that by approximating the nonlinear terms in \eqref{eqn:mag08} and \eqref{eqn:mag09}, we in no way imply that the ratio of voltages between different phases in the resulting linearized model are constant (i.e. the linearized model captures unbalanced voltages). 

Applying these approximations for $\mathbf{H}_{mn}$ and $\Gamma_{n}$ to \eqref{eqn:mag09}, we arrive at a linear system of equations:
\begin{equation}
	\mathbf{E}_{m} \approx \mathbf{E}_{n} + 2 \mathbf{M}_{mn} \mathbf{P}_{mn} - 2 \mathbf{N}_{mn} \mathbf{Q}_{mn} \q ,
    \label{eqn:mag10}
\end{equation}
\begin{equation}
	\mathbf{M}_{mn} = \Re \left\{ A \circ \mathbf{Z}_{mn}^{*} \right\}, \q \mathbf{N}_{mn} = \Im \left\{ A \circ \mathbf{Z}_{mn}^{*} \right\} \q .
    \label{eqn:MmnNmn}
\end{equation}

\noindent The matrices $\mathbf{M}_{mn}$ and $\mathbf{N}_{mn}$ are modified impedance matrices, where the off-diagonal elements are rotated by $\pm 120 \degree$ (see \eqref{eqn:gammaalpha}). The diagonal entries of $\mathbf{M}_{mn}$ are $r_{mn}^{\phi \phi}$. Off-diagonal entries of $\mathbf{M}_{mn}$ are $\frac{1}{2} \left( -r_{mn}^{\phi \psi} + \sqrt{3} x_{mn}^{\phi \psi} \right)$ for $(\phi, \psi) \in \left\{ab, bc, ca \right\}$, and $\frac{1}{2} \left( -r_{mn}^{\phi \psi} - \sqrt{3} x_{mn}^{\phi \psi} \right)$ for $(\phi, \psi) \in \left\{ ac, ba, cb \right\}$. Diagonal entries of $\mathbf{N}_{mn}$ are $-x_{mn}^{\phi \phi}$. Off-diagonal entries of $\mathbf{N}_{mn}$ are $\frac{1}{2} \left( x_{mn}^{\phi \psi} + \sqrt{3} r_{mn}^{\phi \psi} \right)$ for $(\phi, \psi) \in \left\{ ab, bc, ca \right\}$, and $\frac{1}{2} \left( x_{mn}^{\phi \psi} - \sqrt{3} r_{mn}^{\phi \psi} \right)$ for $(\phi, \psi) \in \left\{ ac, ba, cb \right\}$.


\subsection{Motivating Intermezzo}
\label{sbsec:motivation}

Now that we have set up the necessary equations, it is possible to further demonstrate the need of controlling phase angle for switching operations. Consider the following illustrative example, where we assume a single phase line and therefore omit superscripts denoting phase. For two nodes $m$ and $n$ that are not connected by a line, no current flows between the nodes. However should a switch between the two nodes be closed, the power at node $n$ is $S_{mn} = V_{n} \left( V_{m} - V_{n} \right)^{*} Y_{mn}^{*}$, where $Y_{mn} = g_{mn} + j b_{mn}$ is the admittance of the line. Assuming the voltage magnitudes at nodes $m$ and $n$ are equal, the line power at node $n$ can be written as:
\begin{equation*}
\begin{aligned}
	S_{mn} = & \left| V_{n} \right|^{2} \left( g_{mn} ( \cos ( \theta_{mn} ) - 1 ) - b_{mn} \sin ( \theta_{mn} ) \right) \ldots \\
    & + j \left| V_{n} \right|^{2} \left( b_{mn} (1 -  \cos ( \theta_{mn} ) ) - g_{mn} \sin ( \theta_{mn} ) \right) \q ,
\end{aligned}
\end{equation*}

\noindent where $\theta_{mn} = \theta_{m} - \theta_{n}$. It is clear that even with equal voltage magnitudes at nodes $m$ and $n$, larger voltage angle differences will cause increased real and reactive power flows. This highlights the importance of the ability to control voltage angle, which present OPF formulations lack.

\subsection{Voltage Phase Angle Equations}
\label{sbsec:angle}

We now derive an extension of the power and voltage magnitude system that relates differences in voltage angles between adjacent nodes to complex power flows.
This derivation builds heavily upon the analysis of Section \ref{sbsec:magnitude}.

The derivation presented here represents a three phase line, $\mathcal{P}_{mn} = \{ a, b, c \}$.
For lines with less than three phases ($\left| \mathcal{P}_{mn} \right| \le 2$), all equations should be indexed by $\mathcal{P}_{mn}$ as \eqref{eqn:KVLmn03} is.

We begin with the Hadamard Product of $\mathbf{V}_{n}$ and the complex conjugate of \eqref{eqn:KVLmn02}:
\begin{equation}
\mathbf{V}_{m}^{*} \circ \mathbf{V}_{n} = \mathbf{V}_{n}^{*} \circ \mathbf{V}_{n} + \left( \mathbf{Z}_{mn} \mathbf{I}_{mn} \right)^{*} \circ \mathbf{V}_{n} \q .
\label{eqn:angle01}
\end{equation}

\noindent From the analysis in Section \ref{sbsec:magnitude}, we substitute both terms on the RHS, and expand the LHS with the polar representations of voltage phasors:
\begin{equation}
\begin{bmatrix}
\left| V_{m}^{a} \right| \left| V_{n}^{a} \right| \angle \left( -\theta_{m}^{a} + \theta_{n}^{a} \right) \\
\left| V_{m}^{b} \right| \left| V_{n}^{b} \right| \angle \left( -\theta_{m}^{b} + \theta_{n}^{b} \right) \\
\left| V_{m}^{c} \right| \left| V_{n}^{c} \right| \angle \left( -\theta_{m}^{c} + \theta_{n}^{c} \right) \\
\end{bmatrix}
= \mathbf{E}_{n} + \left( \Gamma_{n} \circ \mathbf{Z}_{mn}^{*} \right) \mathbf{S}_{mn} \q .
\label{eqn:angle02}
\end{equation}

\noindent We negate \eqref{eqn:angle02} and take the imaginary component of both sides:
\begin{equation}
\begin{aligned}
\begin{bmatrix}
\left| V_{m}^{a} \right| \left| V_{n}^{a} \right| \sin \left( \theta_{m}^{a} - \theta_{n}^{a} \right) \\
\left| V_{m}^{b} \right| \left| V_{n}^{b} \right| \sin \left( \theta_{m}^{b} - \theta_{n}^{b} \right) \\
\left| V_{m}^{c} \right| \left| V_{n}^{c} \right| \sin \left( \theta_{m}^{c} - \theta_{n}^{c} \right) \\
\end{bmatrix}
& = -\Im \left\{ \left( \Gamma_{n} \circ \mathbf{Z}_{mn}^{*} \right) \mathbf{S}_{mn} \right\} \\
\ldots & = - \mathbf{N}_{mn} \mathbf{P}_{mn} - \mathbf{M}_{mn} \mathbf{Q}_{mn} \q .
\end{aligned}
\label{eqn:angle03}
\end{equation}

\noindent where $\mathbf{M}_{mn}$ and $\mathbf{N}_{mn}$ are defined as in \eqref{eqn:mag09}.

Inspection of the voltage angle equation reveals some interesting similarities compared to the voltage magnitude equations \eqref{eqn:mag09}.
The RHS of \eqref{eqn:mag09} and \eqref{eqn:angle03} are the real and imaginary parts of the same argument (except for a scaling factor of one-half).

To simplify \eqref{eqn:angle03}, we apply the same assumptions of \cite{gan2014convex} to the RHS of \eqref{eqn:angle03}.
Second we assume that $\theta_{m}^{\phi} - \theta_{n}^{\phi}$ is sufficiently small such that the small angle approximation holds, so that $ \sin \left( \theta_{m}^{\phi} - \theta_{n}^{\phi} \right) \approx \theta_{m}^{\phi} - \theta_{n}^{\phi} \qq \forall \phi \in \mathcal{P}_{mn}, \qq \forall (m,n) \in \mathcal{E}$.
Lastly, we approximate the product of voltage magnitudes on the LHS of \eqref{eqn:angle03} as a constant (more specifically, unity), so that $\left| V_{m}^{\phi} \right| = \left| V_{n}^{\phi} \right| = 1 \qq \forall \phi \in \mathcal{P}_{mn}, \qq \forall (m,n) \in \mathcal{E}$.
With these three assumptions applied to \eqref{eqn:angle03}, we arrive at:
\begin{equation}
\Theta_{m} \approx \Theta_{n} - \mathbf{N}_{mn} \mathbf{P}_{mn} - \mathbf{M}_{mn} \mathbf{Q}_{mn} \q ,
\label{eqn:angle04}
\end{equation}

\noindent with $\Theta_{m} = \left[ \theta_{m}^{a} , \qq \theta_{m}^{b} , \qq \theta_{m}^{c} \right]^{T}$, and $\mathbf{M}_{mn}$ and $\textbf{N}_{mn}$ defined by \eqref{eqn:MmnNmn}.

The accuracy of these approximations in modeling system power flows and voltages will be explored in Section \ref{sec:montecarlo}.

\subsection{Linearized Unbalanced Power Flow Model}
\label{sbsec:linearmodel}

We now present the full set of equations that comprise a linearized model for unbalanced power flow. Equations for lines $(m,n) \in \mathcal{E}$, \eqref{eqn:maglin} - \eqref{eqn:MmnNmnlin}, should be indexed by line phases $\mathcal{P}_{mn}$ as in \eqref{eqn:KVLmn03}.

\begin{tcolorbox}[ams gather]
\text{Per phase node complex load} \nonumber
\\
\forall \phi \in \mathcal{P}_{m}, \q \forall m \in \mathcal{N} \nonumber
\\
s_{m}^{\phi} \left( V_{m}^{\phi} \right) = \left( \beta_{S,m}^{\phi} + \beta_{Z,m}^{\phi} E_{m}^{\phi} \right) d_{m}^{\phi} + w_{m}^{\phi} - j c_{m}^{\phi}
\label{eqn:sVlin01}
\\
\text{Node power flow} \nonumber
\\
\forall m \in \mathcal{N} \nonumber
\\
\sum_{l:(l,m) \in \mathcal{E}} \mathbf{S}_{lm} \approx \mathbf{s}_{m} + \sum_{n:(m,n) \in \mathcal{E}} \mathbf{S}_{mn}
\label{eqn:powerlin}
\\
\text{Magnitude and angle equations for lines} \nonumber
\\
\forall (m,n) \in \mathcal{E} \nonumber
\\
\left[ \mathbf{E}_{m} \approx \mathbf{E}_{n} + 2 \mathbf{M}_{mn} \mathbf{P}_{mn} - 2 \mathbf{N}_{mn} \mathbf{Q}_{mn} \right]_{\mathcal{P}_{mn}}
\label{eqn:maglin}
\\
\left[ \Theta_{m} \approx \Theta_{n} - \mathbf{N}_{mn} \mathbf{P}_{mn} - \mathbf{M}_{mn} \mathbf{Q}_{mn} \right]_{\mathcal{P}_{mn}}
\label{eqn:anglelin}
\\
\mathbf{M}_{mn}= \Re \left\{ A \circ \mathbf{Z}_{mn}^{*} \right\} , \q
\mathbf{N}_{mn}= \Im \left\{ A \circ \mathbf{Z}_{mn}^{*} \right\}
\label{eqn:MmnNmnlin}
\end{tcolorbox}

The accuracy of the approximations in the power and voltage magnitude equations has been investigated in \cite{gan2014convex} and \cite{robbins2016optimal}. In the next section, we perform a Monte Carlo analysis to explore the level of error introduced by the voltage angle equation assumptions.

Both the magnitude equation and angle equation can be extended to lines or switches with transformers. Transformer models are often linear voltage and current, often taking the form $\mathbf{V}_{H} = A \mathbf{V}_{L} + B \mathbf{I}_{L}$ where the subscripts $H$ and $L$ refer to either side of the transformer, and $A$, $B$ are constant matrices \cite{kersting2001distribution}. Modified versions of \eqref{eqn:maglin}, \eqref{eqn:anglelin}, and \eqref{eqn:MmnNmnlin} can be derived, which will be explored in subsequent works.

\section{Model Accuracy Analysis}
\label{sec:montecarlo}

\begin{figure}[t]
\centering
\resizebox {\linewidth} {!} { \begin{tikzpicture}


\node[draw, ultra thick, inner sep=3pt, minimum width=50pt] at (0,0) (inf) {\Large $\infty$};

\path (inf) ++(0,-1.854) node[draw, ultra thick, inner sep=3pt] (650) {\large 650};
\draw[ultra thick] (inf) -- (650) node[pos = 0.5, fill=white] () {\large$abc$};


\path (650) ++(0,-1.854) node[draw, ultra thick, inner sep=3pt] (632) {\large 632};
\draw[ultra thick] (650) -- (632) node[pos = 0.5, fill=white] () {\large$abc$};

\path (632) ++(3,0) node[draw, ultra thick, inner sep=3pt] (633) {\large 633};
\draw[ultra thick] (632) -- (633) node[pos = 0.5, fill=white] () {\large$abc$};

\path (633) ++(3,0) node[draw, ultra thick, inner sep=3pt] (634) {\large 634};
\draw[ultra thick] (633) -- (634) node[pos = 0.5, fill=white] () {\large$abc$};

\path (632) ++(-3,0) node[draw, ultra thick, inner sep=3pt] (645) {\large 645};
\draw[ultra thick] (632) -- (645) node[pos = 0.5, fill=white] () {\large$bc$};

\path (645) ++(-3,0) node[draw, ultra thick, inner sep=3pt] (646) {\large 646};
\draw[ultra thick] (645) -- (646) node[pos = 0.5, fill=white] () {\large$bc$};

\path (632) ++(0,-1.854) node[draw, ultra thick, inner sep=3pt] (671) {\large 671};
\draw[ultra thick] (632) -- (671) node[pos = 0.5, fill=white] () {\large$abc$};

\path (671) ++(3,0) node[draw, ultra thick, inner sep=3pt] (692) {\large 692};
\draw[ultra thick] (671) -- (692) node[pos = 0.5, fill=white] () {\large$abc$};

\path (692) ++(3,0) node[draw, ultra thick, inner sep=3pt] (675) {\large 675};
\draw[ultra thick] (692) -- (675) node[pos = 0.5, fill=white] () {\large$abc$};

\path (671) ++(0,-1.854) node[draw, ultra thick, inner sep=3pt] (680) {\large 680};
\draw[ultra thick] (671) -- (680) node[pos = 0.5, fill=white] () {\large$abc$};

\path (671) ++(-3,0) node[draw, ultra thick, inner sep=3pt] (684) {\large 684};
\draw[ultra thick] (671) -- (684) node[pos = 0.5, fill=white] () {\large$ac$};

\path (684) ++(0,-1.854) node[draw, ultra thick, inner sep=3pt] (652) {\large 652};
\draw[ultra thick] (684) -- (652) node[pos = 0.5, fill=white] () {\large$c$};

\path (684) ++(-3,0) node[draw, ultra thick, inner sep=3pt] (611) {\large 611};
\draw[ultra thick] (684) -- (611) node[pos = 0.5, fill=white] () {\large$a$};




\end{tikzpicture} }
\caption{Network topology of IEEE 13 node feeder with line phases 
shown.}
\label{fig:ieee13node}
\end{figure}
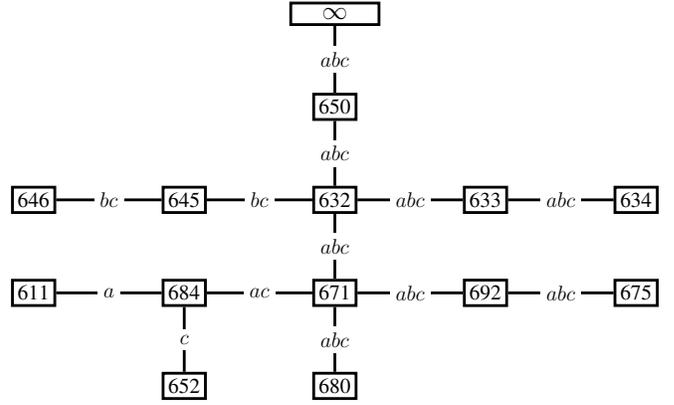

To investigate the accuracy of the approximations derived in the previous section, we perform a Monte Carlo simulation.
We compare the solutions between the non-approximated power flow model, of \eqref{eqn:KVLmn03}, \eqref{eqn:KCLm02}, and \eqref{eqn:sV}, with the linearized power flow model, of \eqref{eqn:sVlin01} - \eqref{eqn:MmnNmnlin}, to analyze and quantify error introduced by the linearizing assumptions on an unbalanced network.
The Monte Carlo simulation was performed on a modified version of the IEEE 13 node test feeder \cite{IEEEtestfeeder}, as seen in Fig.
\ref{fig:ieee13node}.
The voltage regulator between nodes 650 and 632 was omitted, the transformer between nodes 633 and 634 was replaced by a line of configuration 601 (according to \cite{IEEEtestfeeder}, page 5) and length of 50 feet, and the switch between node 671 and 692 was replaced by a line with configuration 601 and length of 50 feet.
The transmission line, denoted in Fig.
\ref{fig:ieee13node} as $\infty$ had a fixed voltage reference of $\mathbf{V}_{\infty} = {\left[ 1 , \q 1 \angle 240 \degree , \q 1 \angle 120 \degree \right]}^{T}$ p.u.

To ensure heterogeneity between loads at different nodes in the network, we vary loading conditions through the following method.
We define $\overline{dr} = \{ 0, 0.01, 0.02, \ldots, 0.15 \}$ as the vector of maximum real component of node load, and $\overline{di} = \{ 0, 0.01, 0.02, \ldots, 0.15 \}$ as the vector of maximum reactive component of node load.
We consider the $\left| \overline{dr} \right| \times \left| \overline{di} \right|$ possible combinations of maximum nodal loading, and at each combination solve power flow for 100 separate scenarios.
For each scenario, real and reactive components of loads were drawn from two separate uniform distributions parametrized by the elements of $\overline{dr}$ and $\overline{di}$, where $d_{n}^{\phi} = \mathcal{U}(0,\overline{dr}) + j \mathcal{U}(0,\overline{di})$.
Loads were then assigned to network locations with spot loads in the original IEEE 13 node feeder documentation (including all single and double phase loads).
All node power demands were assumed to have constant power and constant impedance load parameters of $\beta_{S,n}^{\phi} = 0.85$ and $\beta_{Z,n}^{\phi} = 0.15$.
Power injections from controllable DER were not considered in this experiment, and were therefore set to zero, (i.e. $w_{n}^{\phi}= 0$).

For each scenario (of which there are $15 \times 15 \times 100$), we considered two cases of solving power flow.
The first is non-approximated power flow, henceforth referred to as ``base", composed of Eqs \eqref{eqn:KVLmn03}, \eqref{eqn:KCLm02}, and \eqref{eqn:sV}.
The second is approximated power flow, referred to as ``approximate", comprised of the linearized model \eqref{eqn:sVlin01} - \eqref{eqn:MmnNmnlin}.
Power flow for both the ``base'' and ``approximate'' cases were solved with a Newton-Raphson method adapted from \cite{wasley1974newton}.
The error between base and approximate power flow results are given by \eqref{eqn:magerror} - \eqref{eqn:powerror}, where the $\tilde{\cdot}$ notation indicates the approximate power flow solution.
We define the error in voltage magnitude by \eqref{eqn:magerror}, the error in voltage angle by \eqref{eqn:angerror}, and the error in complex power by \eqref{eqn:powerror}.
These equations capture the maximum absolute (i.e. worst case) error between the base and approximate power flow value across the entire network for a single scenario.
\begin{align}
\varepsilon_{mag} & = \max_{\phi \in \mathcal{P}_{n}, \q n \in \mathcal{N}} \left| \left| V_{n}^{\phi} \right| - \left| \tilde{V}_{n}^{\phi} \right| \right| \label{eqn:magerror} \\
\varepsilon_{angle} & = \max_{\phi \in \mathcal{P}_{n}, \q n \in \mathcal{N}} \left| \angle V_{n}^{\phi} - \angle \tilde{V}_{n}^{\phi} \right| \label{eqn:angerror} \\
\varepsilon_{power} & = \max_{\phi \in \mathcal{P}_{mn}, \q (m,n) \in \mathcal{E}} \left| S_{mn}^{\phi} - \tilde{S}_{mn}^{\phi} \right| \label{eqn:powerror}
\end{align}

In this experiment, substation power, $S_{sub}$, is defined as the sum of the apparent power magnitude delivered by the transmission system in each phase, as in:
\begin{equation}
	S_{sub} = \sum_{\phi \in \{ a, b, c \}} \left| S_{\infty, 650}^{\phi} \right| \q .
\end{equation}

\begin{figure}[t]
\centering
\includegraphics[width=1.0\linewidth]{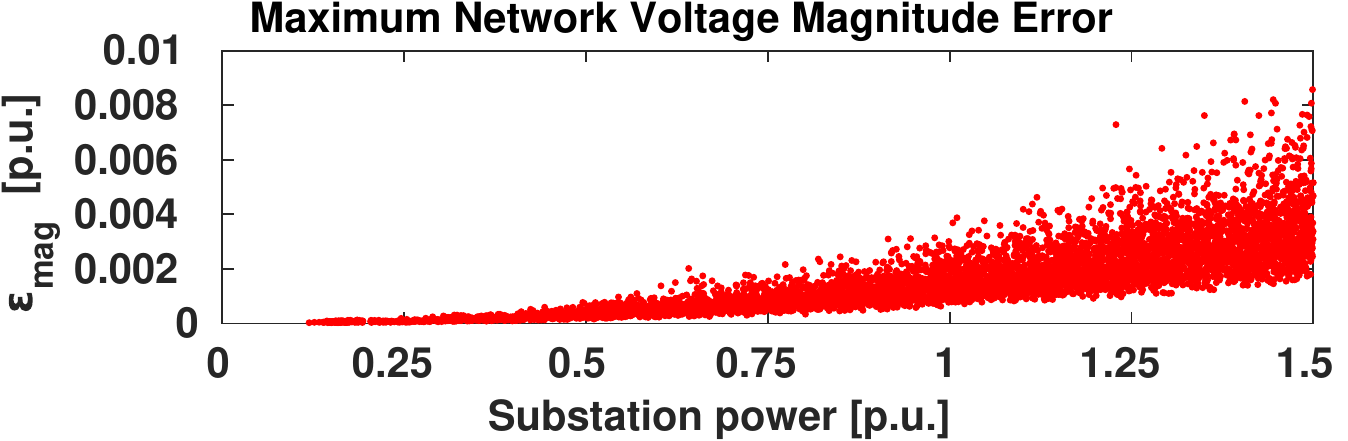}
\caption{Maximum network voltage magnitude error from Monte Carlo simulation, as defined by \eqref{eqn:magerror}}
\label{fig:montecarlomag}
\includegraphics[width=1.0\linewidth]{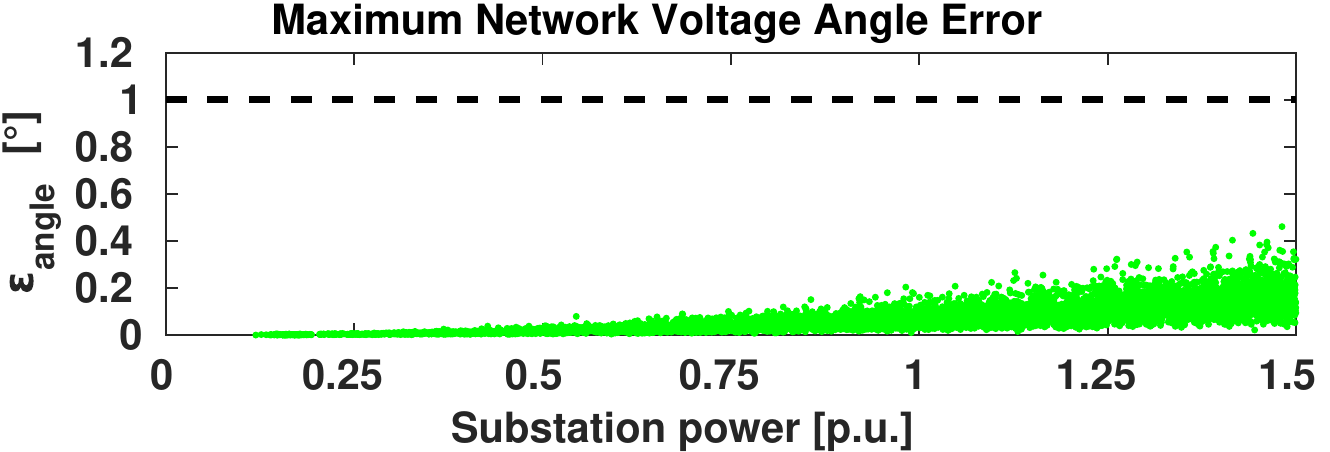}
\caption{Maximum network voltage angle error from Monte Carlo simulation, as defined by \eqref{eqn:angerror}.
The dashed line represents the resolution of a typical synchro-check relay \cite{abbspau140c}.}
\label{fig:montecarloangle}
\includegraphics[width=1.0\linewidth]{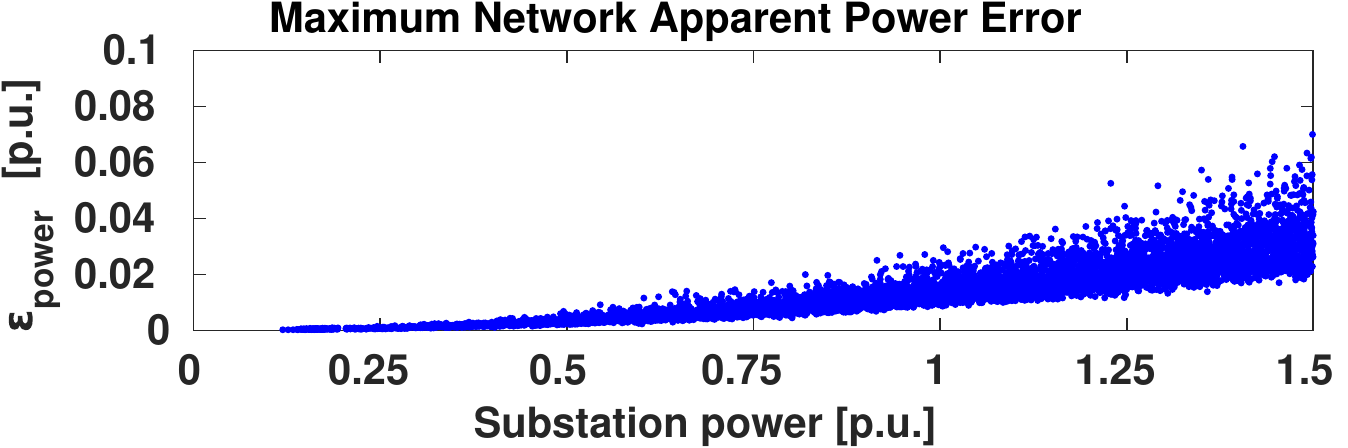}
\caption{Maximum network apparent power error from Monte Carlo simulation, as defined by \eqref{eqn:powerror}.}
\label{fig:montecarlopower}
\end{figure}

Fig.~\ref{fig:montecarlomag} shows voltage magnitude error \eqref{eqn:magerror} for increasing total substation power.
It can clearly be seen that under normal operating conditions ($\le$ 1 p.u. substation rated power) the voltage magnitude error is less than 0.5\%.
For substation loading 1.5 times the rated power, the error remains under 1\%.
Fig.~\ref{fig:montecarloangle} gives the voltage angle error \eqref{eqn:angerror} for increasing total substation power.
At substation rated power, the error typically remains under $0.25 \degree$, much less than the $1 \degree$ resolution of a typical synchro-check relay \cite{abbspau140c}.
Fig.~\ref{fig:montecarlopower} shows the line apparent power error \eqref{eqn:powerror} as a function of substation power.
At substation rated power, we see an error of 2\%.
Although it is observed the error in voltage magnitude, voltage angle, and substation power monotonically increase with substation power, we note that the maximum errors are relatively small under normal operating conditions ($\sim$1 p.u. substation rated power).
%

The results of the preceding analysis reveal that the linearized model becomes increasingly inaccurate has system loading is increased.
However, we note that even under these unlikely conditions, these errors are relatively small and reflect "worst case" conditions.

\section{Phasor Tracking for Switching Operations}
\label{sec:phasortracking}

We now present results of two experiments in which the linearized unbalanced power flow model, \eqref{eqn:sVlin01} - \eqref{eqn:MmnNmnlin}, is incorporated into OPFs with the objective of minimizing the phasor difference between nodes at either end of one or more open switches (we will refer to this as phasor tracking), while regulating system voltage magnitudes to within acceptable limits. The OPF decision variables were DER real and reactive power injections at select nodes, which were capacity constrained. This four-quadrant model of DER was chosen for this experiment in order to analyze the effect of devices that can both source and sink both real and reactive power.

The scenario we discuss is the reconfiguration of distribution networks through switching actions. We do not consider economic activity or optimization in these experiments, rather we assume that DER dispatch for economic purposes is suspended in order to devote said resources to reconfiguration efforts.  In the future, it is feasible to consider re-purposing privately owned or operated assets to protect critical infrastructure in an extreme weather event, cyber attack, or component failure.



\subsection{Previous work on SDP OPFs}
Prior to discussing the OPF, we first wish to include a note regarding the formulation of the phasor tracking problem as an SDP. In an effort to compare the result of our approach with an optimal effort (i.e. solving an OPF that uses the exact power flow equations), we investigated extending the work of \cite{dall2013distributed} to optimize voltage phasor differences.


Initial simulations were conducted on simple balanced radial networks with 6 nodes. We found that some of our simulations successfully returned a rank-one solution. However, this success was not repeatable under small changes of simulation parameters, on larger balanced, or any size unbalanced networks.

Our exploration does not conclude that a rank-one SDP (SDR) solution cannot be obtained when voltage phasor differences in unbalanced networks. However, we wished to provide a short note on our preliminary efforts to incorporate \eqref{eqn:angle03} into an OPF, which motivated the development of the linearized unbalanced power flow model in this work.

\subsection{IEEE 13 Node Feeder Test Case}

\begin{figure}[t]
\centering
\resizebox {1.0\linewidth} {!} { \begin{tikzpicture}


\node[draw, ultra thick, inner sep=3pt, minimum width=50pt] at (0,0) (inf) {\LARGE $\infty$};


\path (inf) ++(-5.25,-2) node[draw, ultra thick, inner sep=3pt] (1650) {\Large 1650};
\draw[ultra thick] (inf) -- ++(0,-1) -| (1650);

\path (1650) ++(0,-2) node[draw, ultra thick, inner sep=3pt, red] (1632) {\Large 1632};
\draw[ultra thick] (1650) -- (1632);

\path (1632) ++(2,0) node[draw, ultra thick, inner sep=3pt] (1633) {\Large 1633};
\draw[ultra thick] (1632) -- (1633);

\path (1633) ++(2,0) node[draw, ultra thick, inner sep=3pt] (1634) {\Large 1634};
\draw[ultra thick] (1633) -- (1634);

\path (1632) ++(-2,0) node[draw, ultra thick, inner sep=3pt] (1645) {\Large 1645};
\draw[ultra thick] (1632) -- (1645);

\path (1645) ++(-2,0) node[draw, ultra thick, inner sep=3pt] (1646) {\Large 1646};
\draw[ultra thick] (1645) -- (1646);

\path (1632) ++(0,-2) node[draw, ultra thick, inner sep=3pt] (1671) {\Large 1671};
\draw[ultra thick] (1632) -- (1671);

\path (1671) ++(2,0) node[draw, ultra thick, inner sep=3pt] (1692) {\Large 1692};
\draw[ultra thick] (1671) -- (1692);

\path (1692) ++(2,0) node[draw, ultra thick, inner sep=3pt, red] (1675) {\Large 1675};
\draw[ultra thick] (1692) -- (1675);

\path (1671) ++(0,-2) node[draw, ultra thick, inner sep=3pt] (1680) {\Large 1680};
\draw[ultra thick] (1671) -- (1680);

\path (1671) ++(-2,0) node[draw, ultra thick, inner sep=3pt, red] (1684) {\Large 1684};
\draw[ultra thick] (1671) -- (1684);

\path (1684) ++(0,-2) node[draw, ultra thick, inner sep=3pt] (1652) {\Large 1652};
\draw[ultra thick] (1684) -- (1652);

\path (1684) ++(-2,0) node[draw, ultra thick, inner sep=3pt] (1611) {\Large 1611};
\draw[ultra thick] (1684) -- (1611);

\draw[ultra thick, dashed, red] (1650) ++(-5,0.75) -- ++(10,0) -- ++(0,-6) -- ++(-2,-2) -- ++(-8,0) -- ++(0,8) node[below right] (N1) {\Large $\mathcal{T}_{1}$};




\path (inf) ++(5.25,-2) node[draw, ultra thick, inner sep=3pt] (2650) {\Large 2650};
\draw[ultra thick] (inf) -- ++(0,-1) -| (2650);

\path (2650) ++(0,-2) node[draw, ultra thick, inner sep=3pt, blue] (2632) {\Large 2632};
\draw[ultra thick] (2650) -- (2632);

\path (2632) ++(2,0) node[draw, ultra thick, inner sep=3pt] (2633) {\Large 2633};
\draw[ultra thick] (2632) -- (2633);

\path (2633) ++(2,0) node[draw, ultra thick, inner sep=3pt] (2634) {\Large 2634};
\draw[ultra thick] (2633) -- (2634);

\path (2632) ++(-2,0) node[draw, ultra thick, inner sep=3pt] (2645) {\Large 2645};
\draw[ultra thick] (2632) -- (2645);

\path (2645) ++(-2,0) node[draw, ultra thick, inner sep=3pt] (2646) {\Large 2646};
\draw[ultra thick] (2645) -- (2646);

\path (2632) ++(0,-2) node[draw, ultra thick, inner sep=3pt, blue] (2671) {\Large 2671};
\draw[ultra thick] (2632) -- (2671);

\path (2671) ++(2,0) node[draw, ultra thick, inner sep=3pt] (2692) {\Large 2692};
\draw[ultra thick] (2671) -- (2692);

\path (2692) ++(2,0) node[draw, ultra thick, inner sep=3pt] (2675) {\Large 2675};
\draw[ultra thick] (2692) -- (2675);

\path (2671) ++(0,-2) node[draw, ultra thick, inner sep=3pt] (2680) {\Large 2680};
\draw[ultra thick] (2671) -- (2680);

\path (2671) ++(-2,0) node[draw, ultra thick, inner sep=3pt] (2684) {\Large 2684};
\draw[ultra thick] (2671) -- (2684);

\path (2684) ++(0,-2) node[draw, ultra thick, inner sep=3pt] (2652) {\Large 2652};
\draw[ultra thick] (2684) -- (2652);

\path (2684) ++(-2,0) node[draw, ultra thick, inner sep=3pt] (2611) {\Large 2611};
\draw[ultra thick] (2684) -- (2611);

\draw[ultra thick, dashed, blue] (2650) ++(-5,0.75) -- ++(10,0) -- ++(0,-8) -- ++(-8,0) -- ++(-2,2) -- ++(0,6) node[below right] (N2) {\Large $\mathcal{T}_{2}$};


\draw[ultra thick, dashed, black!80] (inf) ++(-10.5,0.5) -- ++(21,0) -- ++(0,-10) -- ++(-21,0) -- ++(0,10) node[below right] (N) {\Large $\mathcal{T}$};

\path (1680) ++(0,-1) ++(4,0) ++(2,0.75) node[draw, circle, inner sep=2.5pt] (sw1) {} ++(0,-0.75) node[draw, circle, inner sep=2.5pt] (sw2) {};
\draw[ultra thick] (1680) -- ++(0,-1) -- ++(4,0) -- ++(sw1);
\draw[ultra thick] (2680) -- ++(0,-1) -- (sw2);

\end{tikzpicture} }
\caption{Networks $\mathcal{T}_{1}$ and $\mathcal{T}_{2}$ connected to the same transmission line, with open switch. Nodes with DER resources are highlighted in red for $\mathcal{T}_{1}$ or blue for $\mathcal{T}_{2}$.}
\label{fig:ieee13mesh}
\end{figure}
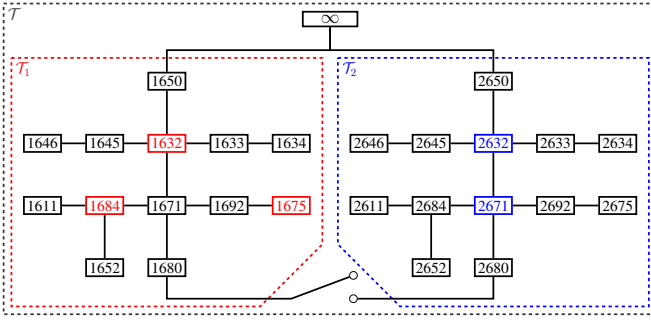


To illustrate the general objective of our phasor tracking strategy, we study an example case of a tie switch connecting two distinct distribution circuits, networks $\mathcal{T}_{1} = (\mathcal{N}_{1},\mathcal{E}_{1})$ and $\mathcal{T}_{2} = (\mathcal{N}_{2},\mathcal{E}_{2})$, as shown in Fig. \ref{fig:ieee13mesh}.
The overall network is represented by $\mathcal{T} = (\mathcal{N},\mathcal{E})$ with $\mathcal{N} = \mathcal{N}_{1} \cup \mathcal{N}_{2} \cup \infty$ and $\mathcal{E} = \mathcal{E}_{1} \cup \mathcal{E}_{2} \cup (\infty,1650) \cup (\infty,2650)$.


Both $\mathcal{T}_{1}$ and $\mathcal{T}_{2}$ were modified versions of the IEEE 13 node test feeder model. Feeder topology, line configuration, line impedance, line length, and spot loads are specified in \cite{IEEEtestfeeder}. To differentiate the elements of the networks between $\mathcal{T}_{1}$ and $\mathcal{T}_{2}$, we prepend the number 1 to the front of nodes within $\mathcal{T}_{1}$ and the number 2 for nodes within $\mathcal{T}_{2}$ (e.g. node 671 of $\mathcal{T}_{1}$ is now 1671 and node 634 of $\mathcal{T}_{2}$ is now 2634). The transmission line was treated as an infinite bus, with a fixed voltage reference of $\mathbf{V}_{\infty} = {\left[ 1 , \qq 1 \angle 240 \degree , \qq 1 \angle 120 \degree \right]}^{T}$ p.u.



The voltage regulators between nodes 1650 and 1632, and between nodes 2650 and 2632, were both omitted. The transformers between nodes 1633 and 1634, and between 2633 and 2634, were both replaced by a line of configuration 601 (according to \cite{IEEEtestfeeder}, page 5) and length of 50 feet. The switches between node 1671 and 1692, and between 2671 and 2692, were both replaced by a line with configuration 601 and length of 50 feet. We placed Wye connected $0.01 + j0.004$ p.u. loads on all phases at node 1680 and 2680.

For both networks, the voltage dependent load model of \eqref{eqn:sV} had parameters $\beta_{S,n}^{\phi} = 0.85$ and $\beta_{Z,n}^{\phi} = 0.15 \qq \forall \phi \in \mathcal{P}_{n}, \qq \forall n \in \mathcal{N}_{1} \cup \mathcal{N}_{2}$. To create a load imbalance between the two networks, we multiplied all loads in $\mathcal{T}_{1}$ by a factor of 0.75, and all loads in $\mathcal{T}_{2}$ by a factor of 1.5. An open switch was placed between node 1680 of $\mathcal{N}_{1}$ and node 2680 of $\mathcal{N}_{2}$, on a line with configuration 601 and length of 500 feet.

Four quadrant capable DER were placed at on all existing phases at nodes $\mathcal{G} = \mathcal{G}_{1} \cup \mathcal{G}_{2}$ with $\mathcal{G}_{1} = \left\{ 1632, 1675, 1684 \right\}$ and $\mathcal{G}_{2} = \left\{ 2632, 2671 \right\}$. We assumed each DER can inject or sink both real and reactive power separately on each phase of the feeder and are only constrained by an apparent power capacity limit on each phase of 0.05 p.u, such that $\overline{w}_{n}^{\phi} = 0.05, \qq \forall \phi \in \mathcal{P}_{n}, \forall n \in \mathcal{G}$ and $w_{n}^{\phi} = \overline{w}_{n}^{\phi} = 0 \qq \forall \phi \in \mathcal{P}_{n}, \forall n \mathcal{N} \setminus \mathcal{G}$. We consider 4 quadrant DER operation to generalize the OPF, such that DER can represent any actionable DER with actionable real power such as PV or wind generation or load vehicle-to-grid operations, load-shedding, or charging/discharging of distributed battery storage.

To capture the effect of the proposed OPF interacting with other devices performing voltage regulation, the following experiment includes other DER performing Volt-VAr control (VVC) that is not controlled by the proposed OPF.  VVC enabled DER was added on all existing phases at the following nodes: 1632, 1680, 2632, 2680. VVC were assumed to have a piecewise linear mapping between node voltage and VAr output (see \cite{inverter2016} for further details on Volt-VAr control):
\begin{equation}
\begin{aligned}
q_{n}^{\phi} & ( V_{n}^{\phi} ) = \ldots \\
& \begin{cases}
\underline{q}, & \left| V_{n}^{\phi} \right| \leq \underline{V} \\
\frac{\overline{q} - \underline{q}}{\overline{V} - \underline{V}} \left( \left| V_{n}^{\phi} \right| - \underline{V} \right) + \underline{q}, & \underline{V} \leq \left| V_{n}^{\phi} \right| \leq \overline{V} \\
\overline{q}, & \left| V_{n}^{\phi} \right| \geq \overline{V}
\end{cases}
\q .
\end{aligned}
\label{eqn:vvc}
\end{equation}

\begin{figure}[t]
\centering
\includegraphics[width=\linewidth]{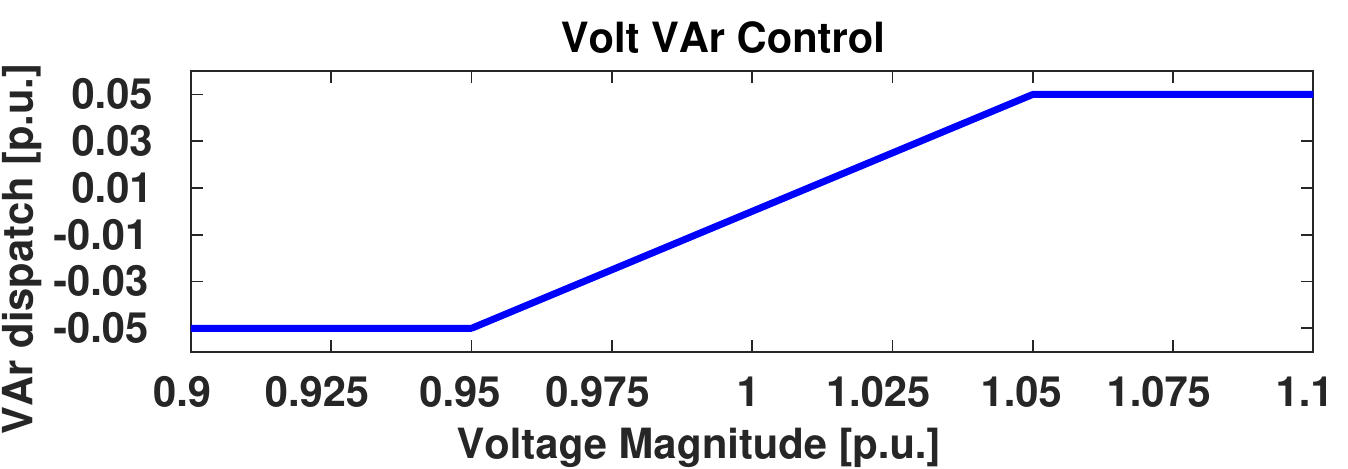}
\caption{Volt-VAr Control.}
\label{fig:vvc}
\end{figure}

All VVC had the following parameters: $\overline{q} = 0.05$ p.u., $\underline{q} = -0.05$ p.u., $\overline{V} = 1.05$ p.u., $\underline{V} = 0.95$ p.u. The VVC model can be seen in Fig. \ref{fig:vvc}. As the linearization used in this work models the square of the voltage magnitude, a Taylor expansion of the form, $|V_{n}^{\phi}| = \sqrt{E_{n}^{\phi}} \approx \frac{1}{2} (1 + E_{n}^{\phi})$, was used to incorporate the piecewise linear VVc model into the OPF as an equality constraint.

In this experiment, our objective was to close an open switch between $\mathcal{T}_{1}$ and $\mathcal{T}_{2}$ thereby connecting the two networks over line $(1680, 2680)$. To minimize large power transfers across the switch upon closing, we desired to match the voltage phasors at the ends of the open switch. To this end, we proposed the following OPF to minimize the voltage phasor difference between one or more nodes and the respective reference at each node, while providing feeder voltage support:
\begin{align}
\begin{aligned}
& \underset{\substack{u_{n}^{\phi}, v_{n}^{\phi}, E_{n}^{\phi} \\ \theta_{n}^{\phi}, P_{n}^{\phi}, Q_{n}^{\phi}}} {\text{minimize}}
& & \rho_{E} C_{E} + \rho_{\theta} C_{\theta} + \rho_{w} C_{w} \\
& \text{subject to}
& & \eqref{eqn:sVlin01}, \eqref{eqn:powerlin}, \eqref{eqn:maglin}, \eqref{eqn:anglelin}, \eqref{eqn:MmnNmnlin} \\
& & & \underline{E} \le E_{n}^{\phi} \le \overline{E} \qq \forall \phi \in \mathcal{P}_{n}, \q \forall n \in \mathcal{N} \\
& & & \mathbf{E}_{\infty} = \left[ 1, \qq 1, \qq 1 \right]^{T} \\
& & & \Theta_{\infty} = \left[ 0, \qq -2\pi/3, \qq 2\pi/3 \right]^{T} \\
& & & \left| w_{n}^{\phi} \right| \le \overline{w}_{n}^{\phi} \qq \forall \phi \in \mathcal{P}_{n}, \q \forall n \in \mathcal{G} \\
&&& q_{n}^{\phi} ( E_{n}^{\phi} ) = \frac{\overline{q} - \underline{q}}{\overline{V} - \underline{V}} \left( \frac{\left( 1 + E_{n}^{\phi} \right)}{2} - \underline{V} \right) + \underline{q}
\q ,
\end{aligned}
\label{eqn:OPF1}
\end{align}
 \noindent where
\begin{align}
C_{E} &= \sum_{\phi \in \mathcal{P}_{k_{1},k_{2}}} {\left( E_{k_{1}}^{\phi} - E_{k_{2}}^{\phi} \right)}^{2}
\q , \label{eqn:OPFCE} \\
C_{\theta} &= \sum_{\phi \in \mathcal{P}_{k_{1},k_{2}}} {\left( \theta_{k_{1}}^{\phi} - \theta_{k_{2}}^{\phi} \right)}^{2}
\q , \label{eqn:OPFCtheta} \\
C_{w} &= \sum_{n \in \mathcal{G}} \sum_{\phi \in \mathcal{P}_{n}} {\left| w_{n}^{\phi} \right|}^{2}
\q . \label{eqn:OPFCw}
\end{align}

The OPF objective function is a weighted sum of three terms, where $k_{1} = 1680$ and $k_{2} = 2680$: $C_{E}$ is the sum of squared voltage magnitude differences squared, $C_{\theta}$ is the sum of voltage angle differences squared, and $C_{w}$ is the sum of the squared magnitudes of all DER dispatch, to avoid applying excessive amounts of control. Constraints of lower and upper voltage magnitude bounds were imposed as $0.95 \le \left| V_{n}^{\phi} \right| \le 1.05 \qq \forall \phi \in \mathcal{P}_{n}, \forall n \in \mathcal{N}$ such that $\underline{E} = 0.9025$ and $\overline{E} = 1.1025$. Additionally, DER dispatch is constrained by its apparent power capacity, $\overline{w}_{n}^{\phi}$.


Results from this experiment can be seen in Table \ref{tab:ieee13switch}. We consider three cases: In the ``No Control'' (NC) case, all DER dispatch is $0$. In the ``Magnitude Control'' (MC), the optimal DER dispatch is obtained solving \eqref{eqn:OPF1} with $\rho_{E} = 1000$, $\rho_{\theta} = 0$, and $\rho_{w} = 1$. In the ``Phasor Control'' (PC) case, the optimal DER dispatch is obtained solving \eqref{eqn:OPF1} with $\rho_{E} = 1000$, $\rho_{\theta} = 1000$, and $\rho_{w} = 1$. We simulate the ``MC" case to illustrate that controlling variables DER real power, DER reactive power, and node voltage magnitude cannot always minimize voltage phasor potential. The objective function weightings $\rho_{E}$, $\rho_{\theta}$, and $\rho_{w}$ were chosen due to normalize effect of system impedances on the voltage magnitude and phase angle components of the objective.  Many impedance terms are on the order of 0.001 p.u.
and, in the linearized model, squared voltage magnitudes and phase angles are essentially linear combinations of power injections scaled by these impedances.

The results of this experiment are given in Table \ref{tab:ieee13switch} that shows the voltage phasor values on either side of the switch, the voltage magnitude and phase angle differences as well as the resulting steady state power flowing in the line after the switch is closed.  Results are shown for all three cases of control (``No Control", ``Magnitude Control", and ``Phasor Control").



It can clearly be seen that with MC, the difference in per phase voltage magnitudes of 1680 and 2680 were minimized, however the voltage angle difference still remained large and on the order of the NC case. Thus the steady state power across the closed switch, while generally less than the NC case, were still large and on the same order as the NC case. With PC, both the voltage magnitude differences and voltage angle differences were minimized, and the real and reactive components of the steady state power across the switch was several orders of magnitude smaller than the NC or MC case.

Table \ref{tab:ieee13der} shows the optimal active and reactive power dispatched from controllable DER for the ``Phasor Control" case of Table \ref{tab:ieee13switch}. Note that these DER dispatch values imply a sub-optimal cost-minimizing configuration. However, phasors at node pairs across a closed switch need not be tracked indefinitely. Rather, phasor matching during a brief interval would suffice to safely perform a switching operation, after which time the dispatch can be changed arbitrarily fast (within the constraints of the new topology). The additional energy cost would then be negligible.

\begin{table*}[ht]
\centering
\caption{SIMULATION RESULTS FOR THREE CONTROL CASES.}
\begin{tabular}{| c | c | c | c | c | c |}
\hline
& Phase $\phi$ & No Control & Magnitude Control & Phasor Control \\
\hline
\multirow{3}{*}{\shortstack[]{Node 1680 Voltage Phasor [p.u.] \\ $V_{1680}^{\phi}$}}
& $a$ & $0.9890 \angle {-1.5997}^{\circ}$ & $0.9789 \angle {-1.8451}^{\circ}$ & $0.9793 \angle {-2.6921}^{\circ}$ \\
& $b$ & $0.9965 \angle {-120.7789}^{\circ}$ & $0.9938 \angle {-120.8093}^{\circ}$ & $0.9943 \angle {-121.0344}^{\circ}$ \\
& $c$ & $0.9825 \angle {118.4644}^{\circ}$ & $0.9656 \angle {118.5646}^{\circ}$ & $0.9649 \angle {117.5846}^{\circ}$ \\
\hline
\multirow{3}{*}{\shortstack[]{Node 2680 Voltage Phasor [p.u.] \\ $V_{2680}^{\phi}$}}
& $a$ & $0.9727 \angle {-3.2141}^{\circ}$ & $0.9780 \angle {-3.1611}^{\circ}$ & $0.9791 \angle {-2.7066}^{\circ}$ \\
& $b$ & $0.9888 \angle {-121.4874}^{\circ}$ & $0.9942 \angle {-121.4237}^{\circ}$ & $0.9944 \angle {-121.0353}^{\circ}$ \\
& $c$ & $0.9552 \angle {117.0104}^{\circ}$ & $0.9642 \angle {117.0255}^{\circ}$ & $0.9642 \angle {117.5883}^{\circ}$ \\
\hline
\multirow{3}{*}{\shortstack[]{Voltage Magnitude Difference [p.u.] \\ $\left| V_{1680}^{\phi} \right| - \left| V_{2680}^{\phi} \right|$}}
& $a$ & 0.0163 & 0.0009 & 0.0002 \\
& $b$ & 0.0077 & -0.0004 & -0.0001 \\
& $c$ & 0.0273 & 0.0014 & 0.0007 \\
\hline
\multirow{ 3}{*}{\shortstack[]{Voltage Angle Difference [${}^{\circ}$] \\ $\theta_{1680}^{\phi} - \theta_{2680}^{\phi}$}}
& $a$ & 1.6144 & 1.3160 & 0.0144 \\
& $b$ & 0.7085 & 0.6144 & 0.0010 \\
& $c$ & 1.4540 & 1.5391 & -0.0038 \\
\hline
\multirow{ 3}{*}{\shortstack[]{Steady State Line Power [p.u.] \\ $S_{1680,2680}^{\phi}$}}
& $a$ & 1.6476 + j0.6027 & 1.2697 - j0.2390 & 0.0051 + j0.0108 \\
& $b$ & 1.1643 + j0.5862 & 0.8582 - j0.1165 & 0.0055 + j0.0021 \\
& $c$ & 1.6205 + j0.6595 & 1.1944 - j0.4384 & 0.0100 + j0.0242 \\
\hline
\end{tabular}
\label{tab:ieee13switch}
\end{table*}

\begin{table}[ht]
\centering
\caption{OPTIMAL PER UNIT DER DISPATCH FROM \eqref{eqn:OPF1}.}
\begin{tabular}{| l | c | c | c |}
\hline
Node & Phase $a$, $w_{n}^{a}$ [p.u.] & Phase $b$, $w_{n}^{b}$ [p.u.] & Phase $c$, $w_{n}^{c}$ [p.u.] \\
\hline
1632 & 0.0324 + j0.0171 & 0.0302 + j0.0195 & 0.0385 + j0.0258 \\
1675 & 0.0448 + j0.0222 & 0.0428 + j0.0259 & 0.0418 + j0.0274 \\
1684 & 0.0448 + j0.0222 & 0 & 0.0418 + j0.0274 \\
2632 & -0.0323 - j0.0169 & -0.0302 - j0.0193 & -0.0384 - j0.0252 \\
2671 & -0.0449 - j0.0220 & -0.0428 - j0.0258 & -0.0422 - j0.0269 \\
\hline
\end{tabular}
\label{tab:ieee13der}
\end{table}

\subsection{Network Simulation with Multiple Switching Actions}


\begin{figure}[t]
\centering
\resizebox {1.0\linewidth} {!} { \begin{tikzpicture}


\node[draw, thick, inner sep=2.5pt, minimum width=50pt] at (0,0) (inf) {\huge $\infty$};

\path (inf) ++(-6.25,-3) node[draw, thick, inner sep=2.5pt] (1799) {\LARGE 1799};
\draw[ultra thick] (inf) -- ++(0,-1.5) -| (1799);

\path (1799) ++(0,-2) node[draw, thick, inner sep=2.5pt] (1701) {\LARGE 1701};
\draw[ultra thick] (1799) -- (1701);

\path (1701) ++(0,-2) node[draw, thick, inner sep=2.5pt, red] (1702) {\LARGE 1702};
\draw[ultra thick] (1701) -- (1702);

\path (1702) ++(-2,0) node[draw, thick, inner sep=2.5pt] (1705) {\LARGE 1705};
\draw[ultra thick] (1702) -- (1705);

\path (1705) ++(0,2) node[draw, thick, inner sep=2.5pt] (1712) {\LARGE 1712};
\draw[ultra thick] (1705) -- (1712);

\path (1705) ++(-2,0) node[draw, thick, inner sep=2.5pt] (1742) {\LARGE 1742};
\draw[ultra thick] (1705) -- (1742);

\path (1702) ++(2,0) node[draw, thick, inner sep=2.5pt] (1713) {\LARGE 1713};
\draw[ultra thick] (1702) -- (1713);

\path (1713) ++(2,0) node[draw, thick, inner sep=2.5pt, red] (1704) {\LARGE 1704};
\draw[ultra thick] (1713) -- (1704);

\path (1704) ++(0,-2) node[draw, thick, inner sep=2.5pt] (1714) {\LARGE 1714};
\draw[ultra thick] (1704) -- (1714);

\path (1714) ++(0,-2) node[draw, thick, inner sep=2.5pt] (1718) {\LARGE 1718};
\draw[ultra thick] (1714) -- (1718);

\path (1704) ++(2,0) node[draw, thick, inner sep=2.5pt] (1720) {\LARGE 1720};
\draw[ultra thick] (1704) -- (1720);

\path (1720) ++(0,-2) node[draw, thick, inner sep=2.5pt] (1706) {\LARGE 1706};
\draw[ultra thick] (1720) -- (1706);

\path (1706) ++(0,-2) node[draw, thick, inner sep=3pt] (1725) {\LARGE 1725};
\draw[ultra thick] (1706) -- (1725);

\path (1720) ++(0,2) node[draw, thick, inner sep=2.5pt] (1707) {\LARGE 1707};
\draw[ultra thick] (1720) -- (1707);

\path (1707) ++(0,2) node[draw, thick, inner sep=2.5pt, red] (1724) {\LARGE 1724};
\draw[ultra thick] (1707) -- (1724);

\path (1707) ++(-2,0) node[draw, thick, inner sep=2.5pt] (1722) {\LARGE 1722};
\draw[ultra thick] (1707) -- (1722);

\path (1702) ++(0,-4) node[draw, thick, inner sep=2.5pt] (1703) {\LARGE 1703};
\draw[ultra thick] (1702) -- (1703);

\path (1703) ++(-2,0) node[draw, thick, inner sep=2.5pt] (1727) {\LARGE 1727};
\draw[ultra thick] (1703) -- (1727);

\path (1727) ++(-2,0) node[draw, thick, inner sep=2.5pt] (1744) {\LARGE 1744};
\draw[ultra thick] (1727) -- (1744);

\path (1744) ++(0,-2) node[draw, thick, inner sep=2.5pt] (1728) {\LARGE 1728};
\draw[ultra thick] (1744) -- (1728);

\path (1744) ++(0,2) node[draw, thick, inner sep=2.5pt, red] (1729) {\LARGE 1729};
\draw[ultra thick] (1744) -- (1729);

\path (1703) ++(0,-2) node[draw, thick, inner sep=2.5pt] (1730) {\LARGE 1730};
\draw[ultra thick] (1703) -- (1730);

\path (1730) ++(0,-2) node[draw, thick, inner sep=2.5pt] (1709) {\LARGE 1709};
\draw[ultra thick] (1730) -- (1709);

\path (1709) ++(2,0) node[draw, thick, inner sep=2.5pt] (1731) {\LARGE 1731};
\draw[ultra thick] (1709) -- (1731);

\path (1709) ++(0,-2) node[draw, thick, inner sep=2.5pt] (1775) {\LARGE 1775};
\draw[ultra thick] (1709) -- (1775);

\path (1709) ++(-2,0) node[draw, thick, inner sep=2.5pt] (1708) {\LARGE 1708};
\draw[ultra thick] (1709) -- (1708);

\path (1708) ++(-2,0) node[draw, thick, inner sep=2.5pt, red] (1732) {\LARGE 1732};
\draw[ultra thick] (1708) -- (1732);

\path (1708) ++(0,-2) node[draw, thick, inner sep=2.5pt] (1733) {\LARGE 1733};
\draw[ultra thick] (1708) -- (1733);

\path (1733) ++(0,-2) node[draw, thick, inner sep=2.5pt] (1734) {\LARGE 1734};
\draw[ultra thick] (1733) -- (1734);

\path (1734) ++(-2,0) node[draw, thick, inner sep=2.5pt] (1710) {\LARGE 1710};
\draw[ultra thick] (1734) -- (1710);

\path (1710) ++(0,-2) node[draw, thick, inner sep=2.5pt, red] (1735) {\LARGE 1735};
\draw[ultra thick] (1710) -- (1735);

\path (1710) ++(0,2) node[draw, thick, inner sep=2.5pt] (1736) {\LARGE 1736};
\draw[ultra thick] (1710) -- (1736);

\path (1734) ++(0,-2) node[draw, thick, inner sep=2.5pt, red] (1737) {\LARGE 1737};
\draw[ultra thick] (1734) -- (1737);

\path (1737) ++(2,0) node[draw, thick, inner sep=2.5pt] (1738) {\LARGE 1738};
\draw[ultra thick] (1737) -- (1738);

\path (1738) ++(2,0) node[draw, thick, inner sep=2.5pt, red] (1711) {\LARGE 1711};
\draw[ultra thick] (1738) -- (1711);

\path (1711) ++(2,0) node[draw, thick, inner sep=2.5pt] (1741) {\LARGE 1741};
\draw[ultra thick] (1711) -- (1741);

\path (1711) ++(0,2) node[draw, thick, inner sep=2.5pt] (1740) {\LARGE 1740};
\draw[ultra thick] (1711) -- (1740);

\draw[ultra thick, dashed, red] (1799) ++(-5,1) -- ++(12,0) -- ++(0,-20) -- ++(-12,0) -- ++(0,20) node[below right] (N1) {\LARGE $\mathcal{T}_{1}$};


\path (inf) ++(6.25,-3) node[draw, thick, inner sep=2.5pt] (2799) {\LARGE 2799};
\draw[ultra thick] (inf) -- ++(0,-1.5) -| (2799);

\path (2799) ++(0,-2) node[draw, thick, inner sep=2.5pt] (2701) {\LARGE 2701};
\draw[ultra thick] (2799) -- (2701);

\path (2701) ++(0,-2) node[draw, thick, inner sep=2.5pt, blue] (2702) {\LARGE 2702};
\draw[ultra thick] (2701) -- (2702);

\path (2702) ++(-2,0) node[draw, thick, inner sep=2.5pt] (2705) {\LARGE 2705};
\draw[ultra thick] (2702) -- (2705);

\path (2705) ++(0,2) node[draw, thick, inner sep=2.5pt] (2712) {\LARGE 2712};
\draw[ultra thick] (2705) -- (2712);

\path (2705) ++(-2,0) node[draw, thick, inner sep=2.5pt] (2742) {\LARGE 2742};
\draw[ultra thick] (2705) -- (2742);

\path (2702) ++(2,0) node[draw, thick, inner sep=2.5pt] (2713) {\LARGE 2713};
\draw[ultra thick] (2702) -- (2713);

\path (2713) ++(2,0) node[draw, thick, inner sep=2.5pt, blue] (2704) {\LARGE 2704};
\draw[ultra thick] (2713) -- (2704);

\path (2704) ++(0,-2) node[draw, thick, inner sep=2.5pt] (2714) {\LARGE 2714};
\draw[ultra thick] (2704) -- (2714);

\path (2714) ++(0,-2) node[draw, thick, inner sep=2.5pt] (2718) {\LARGE 2718};
\draw[ultra thick] (2714) -- (2718);

\path (2704) ++(2,0) node[draw, thick, inner sep=2.5pt] (2720) {\LARGE 2720};
\draw[ultra thick] (2704) -- (2720);

\path (2720) ++(0,-2) node[draw, thick, inner sep=2.5pt] (2706) {\LARGE 2706};
\draw[ultra thick] (2720) -- (2706);

\path (2706) ++(0,-2) node[draw, thick, inner sep=3pt] (2725) {\LARGE 2725};
\draw[ultra thick] (2706) -- (2725);

\path (2720) ++(0,2) node[draw, thick, inner sep=2.5pt] (2707) {\LARGE 2707};
\draw[ultra thick] (2720) -- (2707);

\path (2707) ++(0,2) node[draw, thick, inner sep=2.5pt, blue] (2724) {\LARGE 2724};
\draw[ultra thick] (2707) -- (2724);

\path (2707) ++(-2,0) node[draw, thick, inner sep=2.5pt] (2722) {\LARGE 2722};
\draw[ultra thick] (2707) -- (2722);

\path (2702) ++(0,-4) node[draw, thick, inner sep=2.5pt] (2703) {\LARGE 2703};
\draw[ultra thick] (2702) -- (2703);

\path (2703) ++(-2,0) node[draw, thick, inner sep=2.5pt] (2727) {\LARGE 2727};
\draw[ultra thick] (2703) -- (2727);

\path (2727) ++(-2,0) node[draw, thick, inner sep=2.5pt] (2744) {\LARGE 2744};
\draw[ultra thick] (2727) -- (2744);

\path (2744) ++(0,-2) node[draw, thick, inner sep=2.5pt] (2728) {\LARGE 2728};
\draw[ultra thick] (2744) -- (2728);

\path (2744) ++(0,2) node[draw, thick, inner sep=2.5pt, blue] (2729) {\LARGE 2729};
\draw[ultra thick] (2744) -- (2729);

\path (2703) ++(0,-2) node[draw, thick, inner sep=2.5pt] (2730) {\LARGE 2730};
\draw[ultra thick] (2703) -- (2730);

\path (2730) ++(0,-2) node[draw, thick, inner sep=2.5pt] (2709) {\LARGE 2709};
\draw[ultra thick] (2730) -- (2709);

\path (2709) ++(2,0) node[draw, thick, inner sep=2.5pt] (2731) {\LARGE 2731};
\draw[ultra thick] (2709) -- (2731);

\path (2709) ++(0,-2) node[draw, thick, inner sep=2.5pt] (2775) {\LARGE 2775};
\draw[ultra thick] (2709) -- (2775);

\path (2709) ++(-2,0) node[draw, thick, inner sep=2.5pt] (2708) {\LARGE 2708};
\draw[ultra thick] (2709) -- (2708);

\path (2708) ++(-2,0) node[draw, thick, inner sep=2.5pt] (2732) {\LARGE 2732};
\draw[ultra thick] (2708) -- (2732);

\path (2708) ++(0,-2) node[draw, thick, inner sep=2.5pt] (2733) {\LARGE 2733};
\draw[ultra thick] (2708) -- (2733);

\path (2733) ++(0,-2) node[draw, thick, inner sep=2.5pt] (2734) {\LARGE 2734};
\draw[ultra thick] (2733) -- (2734);

\path (2734) ++(-2,0) node[draw, thick, inner sep=2.5pt] (2710) {\LARGE 2710};
\draw[ultra thick] (2734) -- (2710);

\path (2710) ++(0,-2) node[draw, thick, inner sep=2.5pt, blue] (2735) {\LARGE 2735};
\draw[ultra thick] (2710) -- (2735);

\path (2710) ++(0,2) node[draw, thick, inner sep=2.5pt] (2736) {\LARGE 2736};
\draw[ultra thick] (2710) -- (2736);

\path (2734) ++(0,-2) node[draw, thick, inner sep=2.5pt, blue] (2737) {\LARGE 2737};
\draw[ultra thick] (2734) -- (2737);

\path (2737) ++(2,0) node[draw, thick, inner sep=2.5pt] (2738) {\LARGE 2738};
\draw[ultra thick] (2737) -- (2738);

\path (2738) ++(2,0) node[draw, thick, inner sep=2.5pt, blue] (2711) {\LARGE 2711};
\draw[ultra thick] (2738) -- (2711);

\path (2711) ++(2,0) node[draw, thick, inner sep=2.5pt] (2741) {\LARGE 2741};
\draw[ultra thick] (2711) -- (2741);

\path (2711) ++(0,2) node[draw, thick, inner sep=2.5pt] (2740) {\LARGE 2740};
\draw[ultra thick] (2711) -- (2740);

\draw[ultra thick, dashed, blue] (2799) ++(-5,1) -- ++(12,0) -- ++(0,-20) -- ++(-12,0) -- ++(0,20) node[below right] (N1) {\LARGE $\mathcal{T}_{2}$};


\draw[ultra thick, dashed, black!75] (1731) -- ++(0,-1) -- ++(12.5,0) node[pos=0.25, fill=white] () {\LARGE Switch 1} -- (2731);

\draw[ultra thick, dashed, black!75] (1725) -- ++(0,-1) -- ++(12.5,0) node[pos=0.75, fill=white] () {\LARGE Switch 2} -- (2725);

\draw[ultra thick, dashed, black!75] (inf) ++(-11.5,0.5) -- ++(25,0) -- ++(0,-23) -- ++(-25,0) -- ++(0,23) node[below right] (N) {\LARGE $\mathcal{T}$};

\end{tikzpicture} }
\caption{Networks $\mathcal{T}_{1}$ and $\mathcal{T}_{2}$ connected to the same transmission line, with two open switches between them. Nodes with DER resources are highlighted in red for $\mathcal{T}_{1}$ or blue for $\mathcal{T}_{2}$.}
\label{fig:ieee37mesh}
\end{figure}
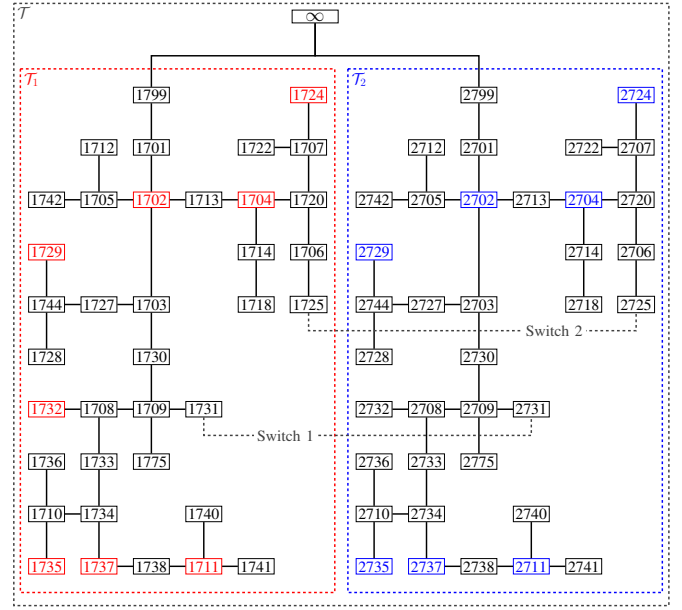


The second experiment considered multiple switching actions on a pair of networks, shown in \ref{fig:ieee37mesh}. Networks $\mathcal{T}_{1}$ and $\mathcal{T}_{2}$, are connected to the same transmission line as shown in Fig. \ref{fig:ieee37mesh}. The overall network is $\mathcal{T} = ( \mathcal{N},\mathcal{E} )$ with $\mathcal{N} = \mathcal{N}_{1} \cup \mathcal{N}_{2} \cup \infty$ and $\mathcal{E} = \mathcal{E}_{1} \cup \mathcal{E}_{2} \cup (\infty,1799) \cup (\infty,2799)$. Both $\mathcal{T}_{1}$ and $\mathcal{T}_{2}$ were modified versions of the IEEE 37 node test feeder model. Feeder topology, line configuration, line impedance, line length, and spot loads are specified in \cite{IEEEtestfeeder}. For clarity, we add the number 1 to the front of nodes within $\mathcal{N}_{1}$ and the number 2 for nodes within $\mathcal{N}_{2}$ (e.g. node 799 of $\mathcal{N}_{1}$ is now 1799 and node 775 of $\mathcal{N}_{2}$ is now 2775). The transmission line was treated as an infinite bus, with a fixed voltage reference of $\mathbf{V}_{\infty} = {\left[ 1 , \qq 1 \angle 240 \degree , \qq 1 \angle 120 \degree \right]}^{T}$ p.u. The transmission line is connected to node $1799$ of $\mathcal{N}_{1}$ (node $799$ in \cite{IEEEtestfeeder}) and $2799$ of $\mathcal{N}_{2}$ (node $799$ in \cite{IEEEtestfeeder}).

The voltage regulators between nodes 1799 and 1701, and between nodes 2799 and 2701, were both omitted. The transformers between nodes 1709 and 1775, and between 1709 and 1775, were both replaced by a line of configuration 724 (according to \cite{IEEEtestfeeder}, page 5) and length of 50 feet. All loads were assumed to be Wye connected on the phase specified in \cite{IEEEtestfeeder}. For both networks, the voltage dependent load model of \eqref{eqn:sV} had parameters $\beta_{S,n}^{\phi} = 0.85$ and $\beta_{Z,n}^{\phi} = 0.15 \qq \forall \phi \in \mathcal{P}_{n}, \qq \forall n \in \mathcal{N}_{1} \cup \mathcal{T}_{2}$. To create a load imbalance between the two networks, we multiplied all loads in $\mathcal{N}_{1}$ by a factor of 1.5, and all loads in $\mathcal{T}_{2}$ by a factor of 1.75.

An open switch was placed between node 1731 of $\mathcal{N}_{1}$ and node 2731 of $\mathcal{N}_{2}$, on a line with configuration 722 and length of 3840 feet. An second open switch was placed between node 1725 of $\mathcal{N}_{1}$ and node 2725 of $\mathcal{N}_{2}$, on a line with configuration 722 and length of 3840 feet.

Four quadrant capable DER were placed at on all existing phases at nodes $\mathcal{G}_{1} = \left\{ 1702, 1704, 1724, 1729, 1732, 1735, 1737, 1711 \right\}$ and $\mathcal{G}_{2} = \left\{ 2702, 2704, 2724, 2729, 2735, 2737, 2711 \right\}$. As in the previous experiment, we assumed each DER can inject or sink both real and reactive power separately on each phase of the feeder and are only constrained by an apparent power capacity limit on each phase of 0.05 p.u, such that $\overline{w}_{n}^{\phi} = 0.05, \qq \forall \phi \in \mathcal{P}_{n}, \forall n \in \mathcal{G}$ and $w_{n}^{\phi} = \overline{w}_{n}^{\phi} = 0 \qq \forall \phi \in \mathcal{P}_{n}, \forall n \notin \mathcal{G}$, where $\mathcal{G} = \mathcal{G}_{1} \cup \mathcal{G}_{2}$.

In this experiment, our objective was to close the two switches sequentially, so as to connect $\mathcal{T}_{1}$ and $\mathcal{T}_{2}$ on two tie lines. Initially both switches were open. DER was dispatched according to \eqref{eqn:OPF1}, with $k_{1} = 1731$ and $k_{2} = 2731$ to minimize the phasor difference between nodes $1731$ and $2731$. After DER was dispatched, the switch was closed and line $(1731,2731)$ was added to $\mathcal{E}$, such that $\mathcal{E} = \mathcal{E} \cup (1731,2731)$.  It is important to note that after the closing of the first switch, the aggregate network is no longer radial.  Another drawback of semidefinite programming approaches is their difficulty in optimizing meshed networks \cite{madani2015convex}. This problem is easily overcome via the use of a linear OPF as discussed herein and demonstrated via controlling the system to enable the next switching action, which is now discussed.

Next, DER was dispatched again according to \eqref{eqn:OPF1}, with $k_{1} = 1725$ and $k_{2} = 2725$ to minimize the phasor difference between nodes $1725$ and $2725$.
The second switch was closed and line $(1725,2725)$ was added to $\mathcal{E}$, such that $\mathcal{E} = \mathcal{E} \cup (1725,2725)$.



For both switching actions, we consider two cases. In the ``No Control'' (NC) Case, all DER dispatch is 0. In the ``Phasor Control'' (PC) Case for the first switching action, the optimal DER dispatch is given by \eqref{eqn:OPF1} with $k_{1} = 1731$ and $k_{2} = 2731$, $\rho_{E} = 1000$, $\rho_{\theta} = 1000$, and $\rho_{w} = 1$. In the ``Phasor Control'' (PC) Case for the second switching action, the optimal DER dispatch is given by \eqref{eqn:OPF1} with $k_{1} = 1725$ and $k_{2} = 2725$, $\rho_{E} = 1000$, $\rho_{\theta} = 1000$, and $\rho_{w} = 1$.

It can clearly be seen in Table \ref{tab:ieee37switch03} that the voltage phasor difference and steady state power flow between nodes 1731 and 2731 is minimized in both switching actions. We note that the steady state power across the closed switch is orders of magnitude smaller than the NC case.

\begin{table*}[ht]
\centering
\caption{SIMULATION RESULTS FOR FIRST AND SECOND SWITCHING ACTIONS.}
\begin{tabular}{| c | c | c | c | c | c | c |}
\hline
& & \multicolumn{2}{|c|}{\shortstack[]{First Switching Action \\ $k_{1} = 1731$, $k_{2} = 2731$}} & \multicolumn{2}{|c|}{\shortstack[]{Second Switching Action \\ $k_{1} = 1725$, $k_{2} = 2725$}} \\
\hline
& Phase $\phi$ & No Control & Phasor Control & No Control & Phasor Control \\
\hline
\multirow{3}{*}{\shortstack[]{Node $k_{1}$ Voltage Phasor [p.u.] \\ $V_{k_{1}}^{\phi}$}}
& $a$ & $0.9784 \angle {-0.5147}^{\circ}$ & $0.9765 \angle {-0.5587}^{\circ}$ & $0.9882 \angle {-0.3003}^{\circ}$ & $0.9872 \angle {-0.3257}^{\circ}$ \\
& $b$ & $0.9909 \angle {-120.3949}^{\circ}$ & $0.9901 \angle {-120.4279}^{\circ}$ & $0.9882 \angle {-120.2812}^{\circ}$ & $0.9833 \angle {-120.3048}^{\circ}$ \\
& $c$ & $0.9749 \angle {119.4164}^{\circ}$ & $0.9727 \angle {119.3669}^{\circ}$ & $0.9882 \angle {119.3006}^{\circ}$ & $0.9787 \angle {119.2415}^{\circ}$ \\
\hline
\multirow{3}{*}{\shortstack[]{Node $k_{2}$ Voltage Phasor [p.u.] \\ $V_{k_{2}}^{\phi}$}}
& $a$ & $0.9747 \angle {-0.6021}^{\circ}$ & $0.9765 \angle {-0.5580}^{\circ}$ & $0.9861 \angle {-0.3508}^{\circ}$ & $0.9871 \angle {-0.3255}^{\circ}$ \\
& $b$ & $0.9893 \angle {-120.4620}^{\circ}$ & $0.9901 \angle {-120.4290}^{\circ}$ & $0.9861 \angle {-120.3289}^{\circ}$ & $0.9832 \angle {-120.3052}^{\circ}$ \\
& $c$ & $0.9706 \angle {119.3172}^{\circ}$ & $0.9727 \angle {119.3669}^{\circ}$ & $0.9861 \angle {119.1815}^{\circ}$ & $0.9786 \angle {119.2407}^{\circ}$ \\
\hline
\multirow{3}{*}{\shortstack[]{Voltage Magnitude Difference [p.u.] \\ $\left| V_{k_{1}}^{\phi} \right| - \left| V_{k_{2}}^{\phi} \right|$}}
& $a$ & 0.0037 & 0.0000 & 0.0020 & 0.0001 \\
& $b$ & 0.0015 & 0.0000 & 0.0026 & 0.0000 \\
& $c$ & 0.0043 & 0.0000 & 0.0033 & 0.0000 \\
\hline
\multirow{ 3}{*}{\shortstack[]{Voltage Angle Difference [${}^{\circ}$] \\ $\theta_{k_{1}}^{\phi} - \theta_{k_{2}}^{\phi}$}}
& $a$ & 0.0874 & -0.0007 & 0.0505 & -0.0002 \\
& $b$ & 0.0671 & 0.0011 & 0.0477 & 0.0004 \\
& $c$ & 0.0991 & -0.0000 & 0.1190 & 0.0009 \\
\hline
\multirow{ 3}{*}{\shortstack[]{Steady State Line Power Power [p.u.] \\ $S_{k_{1},k_{2}}^{\phi}$}}
& $a$ & 0.0852 + j0.0348 & 0.0005 + j0.0005 & 0.0581 + j0.0252 & 0.0013 + j0.0011 \\
& $b$ & 0.0673 + j0.0221 & 0.0010 + j0.0004 & 0.0808 + j0.0316 & 0.0008 + j0.0008 \\
& $c$ & 0.1100 + j0.0442 & 0.0001 + j0.0005 & 0.0867 + j0.0280 & 0.0009 + j0.0006 \\
\hline
\end{tabular}
\label{tab:ieee37switch03}
\end{table*}

\subsection{Discussion on OPFs using Linearized Unbalanced Power Flow Model}

In this work, the implementation of the OPF to regulate voltage phasor differences was accomplished via a centralized entity that computed an optimal decision which was subsequently disseminated to individual controllable DER.  This particular implementation is scalable to manage switching actions for an arbitrary number of networks, provided the existence of network models and the availability of a communications infrastructure to relay sensing and actuation signals.

In situations where multiple switching actions are desired, as presented in the previous section, it should be possible to develop a methodology to determine an optimal switching sequence that would allow one or more of the voltage phasor differences across a switch to be further minimized.  Using the linearized OPF proposed herein, one could enumerate possible switching actions for small networks, but this approach would breakdown as the number of switches and networks being considered increase.  We intend to explore optimal switching sequences in future works.


Although a limitation of the linearized model is the increased inaccuracy of the linearization under high loading conditions, as discussed in Section \ref{sec:montecarlo}, this limitation can be addressed via an iterative approach akin to sequential quadratic programming, in which approximated terms are updated during iterations of the OPF algorithm.


\section{Conclusion}
\label{sec:conclusion}

Optimization of unbalanced power flow is a challenging topic due to its nonlinear and non-convex nature.
While recent works on SDP relaxations \cite{dall2012optimization, dall2013distributed} have made OPF formulations for unbalanced systems possible, these approaches suffer from restrictions on the possible objectives and a high-dimensional geometrical complexity that impedes feasibility and uniqueness of the solutions.

In this paper, we sought to solve a problem that, to our knowledge, cannot be addressed with SDP techniques.
We build upon our previous work \cite{arnold2015optimal} and that of \cite{gan2014convex}, and \cite{robbins2016optimal} to develop an approximate model for distribution power flow that can be incorporated into convex optimal power flow problems, with the intention of enabling better switching in distribution networks.
To do so, in Section \ref{sbsec:angle}, we developed a model that maps complex power flows into voltage angle differences.
This extended model allows the formulation of OPF problems that manage the entire voltage phasor, rather than only voltage magnitude.


In Section \ref{sec:montecarlo}, we investigated the accuracy of the newly-derived linear model, comparing the results of solving power flow with the system physics and the linearization.
We found that under normal operating conditions, the model leads to magnitude errors less than 0.5\%, angle errors of less than 0.25\degree, and substation power errors of 2\% of the network rated power.

We then incorporated the linear model into an OPF to manage DER assets to enable switching in distribution system.
To accomplish this, an OPF was formulated to minimize the voltage phasor difference across an open switch.
Simulation results demonstrate the effectiveness of the OPF in minimizing voltage phasor difference between two disconnected points in a network.

The ability to switch components into and out of distribution feeders with minimal impact on system operation presents many opportunities to reconfigure distribution systems for a variety of purposes.
Moving forward, we intend to investigate two such applications.
First, we plan to study grid reconfiguration in order to better withstand critical grid events (e.g.
weather-related or other types of disasters).
To solve such a problem, we will most likely need to extend our present OPF formulation into a receding horizon controller, that can optimize over a future time window.
Secondly, as ``clean'' switching may also enable distributed microgrids to coalesce and pool resources to provide ancillary services, we intend to extend this OPF formulation to allow for mixed-integer formulations.

We recognize that it may be difficult to solve a centralized OPF of this type in an on-line fashion due to lack of proper network models and a robust communications system.
In previous works exploring OPF approaches to managing DER we have extensively explored alternative approaches to solving centralized OPFs using model-free and low communication optimization techniques \cite{arnold2017model}.
Moving forward, we will attempt to utilize this work to lessen the information and communications requirement associated with the OPF presented in this paper, thereby allowing optimal voltage phasor management in an on-line setting.

\section{Acknowledgements}

We would like to thank Werner van Westering and his colleagues at distribution utility Alliander, The Netherlands for insightful conversations and suggestions regarding the operation of network reconfiguration.



\bibliography{refnumeric}

\begin{thebibliography}{10}
\providecommand{\url}[1]{#1}
\csname url@samestyle\endcsname
\providecommand{\newblock}{\relax}
\providecommand{\bibinfo}[2]{#2}
\providecommand{\BIBentrySTDinterwordspacing}{\spaceskip=0pt\relax}
\providecommand{\BIBentryALTinterwordstretchfactor}{4}
\providecommand{\BIBentryALTinterwordspacing}{\spaceskip=\fontdimen2\font plus
\BIBentryALTinterwordstretchfactor\fontdimen3\font minus
  \fontdimen4\font\relax}
\providecommand{\BIBforeignlanguage}[2]{{%
\expandafter\ifx\csname l@#1\endcsname\relax
\typeout{** WARNING: IEEEtran.bst: No hyphenation pattern has been}%
\typeout{** loaded for the language `#1'. Using the pattern for}%
\typeout{** the default language instead.}%
\else
\language=\csname l@#1\endcsname
\fi
#2}}
\providecommand{\BIBdecl}{\relax}
\BIBdecl

\bibitem{vonMeier2017precision}
A.~von Meier, E.~Stewart, A.~McEachern, M.~Andersen, and L.~Mehrmanesh,
  ``{Precision Micro-Synchrophasors for Distribution Systems: A Summary of
  Applications},'' \emph{IEEE Transactions on Smart Grid}, vol.~8, no.~6, pp.
  2926--2936, Nov. 2017.

\bibitem{synchrophasor2010}
\BIBentryALTinterwordspacing
``Synchrophasor fact sheet,'' accessed, June-2010. [Online]. Available:
  \url{https://cdn.selinc.com/assets/Literature/Product\%20Literature/Flyers/FS_Synchrophasor_BF_20100617.pdf?v=20150408-131001}
\BIBentrySTDinterwordspacing

\bibitem{ochoa2010angle}
L.~F. Ochoa and D.~H. Wilson, ``Angle constraint active management of
  distribution networks with wind power,'' in \emph{Innovative Smart Gird Tech.
  Conf. Europe (ISGT Europe)}.\hskip 1em plus 0.5em minus 0.4em\relax IEEE,
  Oct. 2010, pp. 1--5.

\bibitem{wang2013pmu}
D.~Wang, D.~Wilson, S.~Venkata, and G.~C. Murphy, ``{PMU}-based angle
  constraint active management on 33k{V} distribution network,'' in \emph{22nd
  Int. Conf. Elect. Distribution}.\hskip 1em plus 0.5em minus 0.4em\relax IET,
  2013, pp. 1--4.

\bibitem{grid2015}
\BIBentryALTinterwordspacing
``Grid modernization multi-year program plan,'' accessed Aug.-2015. [Online].
  Available:
  \url{http://energy.gov/sites/prod/files/2016/01/f28/Grid\%20Modernization\%20Multi-Year\%20Program\%20Plan.pdf}
\BIBentrySTDinterwordspacing

\bibitem{quad2015}
\BIBentryALTinterwordspacing
``Quadrennial energy review: First installment,'' accessed Apr.-2015. [Online].
  Available:
  \url{http://energy.gov/epsa/downloads/quadrennial-energy-review-first-installment}
\BIBentrySTDinterwordspacing

\bibitem{heydt2010next}
G.~T. Heydt, ``The next generation of power distribution systems,''
  \emph{{IEEE} Trans. Smart Grid}, vol.~1, no.~3, pp. 225--235, 2010.

\bibitem{zhang2017distributed}
Y.~Zhang, M.~Hong, E.~Dall'Anese, S.~Dhople, and Z.~Xu, ``Distributed
  controllers seeking {AC} optimal power flow solutions using {ADMM},''
  \emph{{IEEE} Trans. Smart Grid}, vol.~PP, no.~99, pp. 1--1, 2017.

\bibitem{xu2017multi}
Y.~Xu, Z.~Y. Dong, R.~Zhang, and D.~J. Hill, ``Multi-timescale coordinated
  voltage/var control of high renewable-penetrated distribution systems,''
  \emph{{IEEE} Trans. Power Syst.}, vol.~32, no.~6, pp. 4398--4408, Nov. 2017.

\bibitem{arnold2017model}
D.~B. Arnold, M.~D. Sankur, M.~Negrete-Pincetic, and D.~Callaway, ``Model-free
  optimal coordination of distributed energy resources for provisioning
  transmission-level services,'' \emph{{IEEE} Trans. Power Syst.}, vol.~33,
  no.~1, pp. 817--829, Jan. 2018.

\bibitem{antoniadou2017distributed}
K.~E. Antoniadou-Plytaria, I.~N. Kouveliotis-Lysikatos, P.~S. Georgilakis, and
  N.~D. Hatziargyriou, ``Distributed and decentralized voltage control of smart
  distribution networks: Models, methods, and future research,'' \emph{{IEEE}
  Trans. Smart Grid}, vol.~8, no.~6, pp. 2999--3008, Nov. 2017.

\bibitem{dall2012optimization}
E.~Dall'Anese, G.~B. Giannakis, and B.~F. Wollenberg, ``Optimization of
  unbalanced power distribution networks via semidefinite relaxation,'' in
  \emph{North Amer. Power Symp. (NAPS)}.\hskip 1em plus 0.5em minus 0.4em\relax
  IEEE, 2012, pp. 1--6.

\bibitem{dall2013distributed}
E.~Dall'Anese, H.~Zhu, and G.~Giannakis, ``Distributed optimal power flow for
  smart microgrids,'' \emph{{IEEE} Trans. Smart Grid}, vol.~4, no.~3, pp.
  1464--1475, Sept. 2013.

\bibitem{lesieutre2011examining}
B.~C. Lesieutre, D.~K. Molzahn, A.~R. Borden, and C.~L. DeMarco, ``Examining
  the limits of the application of semidefinite programming to power flow
  problems,'' in \emph{Communication, Control, and Computing (Allerton), 2011
  49th Annual Allerton Conference on}.\hskip 1em plus 0.5em minus 0.4em\relax
  IEEE, 2011, pp. 1492--1499.

\bibitem{louca2014nondegeneracy}
R.~Louca, P.~Seiler, and E.~Bitar, ``Nondegeneracy and inexactness of
  semidefinite relaxations of optimal power flow,'' \emph{arXiv preprint
  arXiv:1411.4663}, 2014.

\bibitem{madani2015convex}
R.~Madani, S.~Sojoudi, and J.~Lavaei, ``Convex relaxation for optimal power
  flow problem: Mesh networks,'' \emph{{IEEE} Trans. Power Syst.}, vol.~30,
  no.~1, pp. 199--211, Jan. 2015.

\bibitem{gan2014convex}
L.~Gan and S.~H. Low, ``Convex relaxations and linear approximation for optimal
  power flow in multiphase radial networks,'' in \emph{2014 Power Syst.
  Computation Conf. (PSCC)}.\hskip 1em plus 0.5em minus 0.4em\relax IEEE, Aug.,
  pp. 1--9.

\bibitem{robbins2016optimal}
B.~A. Robbins and A.~D. Dom{\'\i}nguez-Garc{\'\i}a, ``Optimal reactive power
  dispatch for voltage regulation in unbalanced distribution systems,''
  \emph{{IEEE} Trans. Power Syst.}, vol.~31, no.~4, pp. 2903--2913, July 2016.

\bibitem{arnold2015optimal}
D.~B. Arnold, M.~D. Sankur, R.~Dobbe, K.~Brady, D.~S. Callaway, and
  A.~Von~Meier, ``Optimal dispatch of reactive power for voltage regulation and
  balancing in unbalanced distribution systems,'' in \emph{2016 IEEE Power
  Energy Soc. Gen. Meeting (PESGM)}, July, pp. 1--5.

\bibitem{baran1989optimal}
M.~E. Baran and F.~F. Wu, ``Optimal sizing of capacitors placed on a radial
  distribution system,'' \emph{{IEEE} Trans. Power Del.}, vol.~4, no.~1, pp.
  735--743, Jan 1989.

\bibitem{kersting2001distribution}
W.~H. Kersting, \emph{Distribution system modeling and analysis}.\hskip 1em
  plus 0.5em minus 0.4em\relax CRC press, 2001.

\bibitem{IEEEtestfeeder}
\BIBentryALTinterwordspacing
``{IEEE} distribution test feeders,'' accessed May-2015. [Online]. Available:
  \url{http://ewh.ieee.org/soc/pes/dsacom/testfeeders/index.html}
\BIBentrySTDinterwordspacing

\bibitem{wasley1974newton}
R.~Wasley and M.~Shlash, ``Newton-{R}aphson algorithm for 3-phase load flow,''
  in \emph{Proc. of the Inst. of Elec. Eng.}, vol. 121, no.~7.\hskip 1em plus
  0.5em minus 0.4em\relax IET, 1974, pp. 630--638.

\bibitem{abbspau140c}
\BIBentryALTinterwordspacing
``{SPAU} 140{C} synchro-check relay, product guide,'' accessed Nov. 2017.
  [Online]. Available:
  \url{https://search-ext.abb.com/library/Download.aspx?DocumentID=1MRS750421-MBG&LanguageCode=en&DocumentPartId=&Action=Launch}
\BIBentrySTDinterwordspacing

\bibitem{inverter2016}
B.~Seal, ``{Common Functions for Smart Inverters, 4th Ed.}'' Electric Power
  Research Institute, Tech. Rep. 3002008217, 2017.

\end{thebibliography}
\bibliographystyle{IEEEtran}

\end{document}